\documentclass[11pt]{amsart}

\usepackage{hyperref}
\hypersetup{nesting=true,debug=true,naturalnames=true}
\usepackage{graphicx,amssymb,upref}
\usepackage[usenames,dvipsnames]{pstricks}
\usepackage{epsfig}
\usepackage{color}
\usepackage[utf8]{inputenc}
\usepackage[T1]{fontenc}
\usepackage{amsmath, amsthm}
\usepackage[english]{babel}
\usepackage{amssymb}
\usepackage{latexsym}
\usepackage{graphicx,amssymb,upref}
\usepackage[all]{xy}
\usepackage{tikz}
\usetikzlibrary{matrix} 
\usepackage{mathtools}
\newcommand\restr[2]{{% we make the whole thing an ordinary symbol
  \left.\kern-\nulldelimiterspace % automatically resize the bar with \right
  #1 % the function
  \vphantom{\big|} % pretend it's a little taller at normal size
  \right|_{#2} % this is the delimiter
  }}

\date{\today}

\newtheorem{theorem}{Theorem}[section]
\newtheorem{lemma}[theorem]{Lemma}

\newtheorem{corollary}[theorem]{Corollary}
\theoremstyle{definition}
\newtheorem{definition}[theorem]{Definition}
\newtheorem{example}[theorem]{Example}

\newtheorem{remark}[theorem]{Remark}

\newtheorem*{theorem*}{Theorem}
\newtheorem*{lemma*}{Lemma}
\newtheorem*{proposition*}{Proposition}
\newtheorem*{definition*}{Definition}
\newtheorem*{remark*}{Remark}

\DeclareMathOperator{\ot}{\otimes}
\DeclareMathOperator{\co}{\circ}

\DeclareMathOperator{\hbr}{\textnormal{\sf{HBrcd}}}

\title[Hopf bracoids]{Hopf bracoids}
\title[Hopf bracoids]{Hopf bracoids}

\begin{document}
	
\maketitle
	
\begin{center}
	{\bf Jos\'e Manuel Fern\'andez Vilaboa$^{1}$, Ram\'on
	Gonz\'{a}lez Rodr\'{\i}guez$^{2}$ and Brais Ramos P\'erez$^{3}$}.
	\end{center}
	
	\vspace{0.4cm}
	\begin{center}
	{\small $^{1}$ [https://orcid.org/0000-0002-5995-7961].}
	\end{center}
	\begin{center}
	{\small  CITMAga, 15782 Santiago de Compostela, Spain.}
	\end{center}
	\begin{center}
	{\small  Universidade de Santiago de Compostela. Departamento de Matem\'aticas,  Facultade de Matem\'aticas, E-15771 Santiago de Compostela, Spain. 
	\\ email: josemanuel.fernandez@usc.es.}
	\end{center}
	\vspace{0.2cm}
	
	\begin{center}
	{\small $^{2}$ [https://orcid.org/0000-0003-3061-6685].}
	\end{center}
	\begin{center}
	{\small  CITMAga, 15782 Santiago de Compostela, Spain.}
	\end{center}
	\begin{center}
	{\small  Universidade de Vigo, Departamento de Matem\'{a}tica Aplicada II,  E. E. Telecomunicaci\'on,
	E-36310  Vigo, Spain.
	\\email: rgon@dma.uvigo.es}
	\end{center}
    \vspace{0.2cm}

    \begin{center}
   	{\small $^{3}$ [https://orcid.org/0009-0006-3912-4483].}
     \end{center}
	\begin{center}
	{\small  CITMAga, 15782 Santiago de Compostela, Spain. \\}
	\end{center}
    \begin{center}
	{\small  Universidade de Santiago de Compostela. Departamento de Matem\'aticas,  Facultade de Matem\'aticas, E-15771 Santiago de Compostela, Spain. 
	\\email: braisramos.perez@usc.es}
	\end{center}
	\vspace{0.2cm}
	
	%\begin{center}
	%{\small $^{d}$ Corresponding author}
	%\end{center}

\begin{abstract}  
The present article is devoted to introduce the notion of Hopf bracoid in the braided monoidal framework as the quantum version of skew bracoids, which have been presented by Martin-Lyons and Paul J. Truman in \cite{MLT}. Taking into account that Hopf braces are particular examples of Hopf bracoids, in this paper we generalize many properties of Hopf braces to the Hopf bracoid setting and  we obtain a categorical isomorphism between certain full subcategories of the Hopf bracoids category and the 1-cocycles category. 
\end{abstract} 

\vspace{0.2cm}


\vspace{0.2 cm}
\section*{Introduction}
Skew braces were introduced by Guarnieri and Vendramin in \cite{GV} as a tool for finding non-degenerate solutions  not necessarily involutive of the Quantum Yang-Baxter Equation. This particular kind of objects are no more than two different groups, $(G,.)$ and $(G,\star)$, which are compatible with each other according to the following expression:
\begin{equation}\label{condcompatSBrace}
g\star(h.t)=(g\star h).g^{-1}.(g\star t),
\end{equation}
for all $g,h,t\in G$, where $g^{-1}$ denotes the inverse of $g$ for the first group structure $(G,.)$. Carrying out a process of linearization on these objects, the notion of Hopf brace given by I. Angiono, C. Galindo and L. Vendramin in \cite{AGV} is obtained: If $(H,\epsilon,\Delta)$ is a coalgebra, a Hopf brace structure over $H$ is a pair of Hopf algebras,
\begin{align*}
H_{1}=(H,1,\cdot,\epsilon,\Delta, S),\quad H_{2}=(H,1_{\circ},\circ,\epsilon,\Delta, T),
\end{align*}
where $S$ and $T$ are the antipodes, satisfying the following compatibility condition between them:
\begin{align}\label{condcompHBrace}
g\circ(h\cdot t)=(g_{1}\circ h)\cdot S(g_{2})\cdot (g_{3}\circ t)
\end{align}
for all $g,h,t\in H$. Hopf braces are also important from the point of view of the mathematical physics because they give rise to solutions of the Braid Equation as it was proven in \cite[Corollary 2.4]{AGV}.

Since the emergence of skew and Hopf braces, many generalizations of these objects have appeared. An example of these can be the concept of Hopf truss introduced by T. Brzezi\'nski in \cite{BRZ}, whose main difference with the Hopf brace structure is that the compatibility equation \eqref{condcompHBrace} is modified using a cocycle $\sigma\colon H\rightarrow H$. Another generalizations can be the structures of semi braces or weak braces (see \cite{CCS} and \cite{CMMS}, respectively). It is worth highlighting that all of the above-mentioned structures are defined over the same underlying object. Thus, what happens if we try to establish a compatibility condition similar to \eqref{condcompatSBrace} between two groups defined over different underlying sets? This motivation was what gave rise to the notion of skew bracoid proposed by Martin-Lyons and Paul J. Truman in \cite{MLT}. A skew bracoid is a 5-tuple $(G,.,N,\star,\odot)$ such that $(G,.)$ and $(N,\star)$ are groups and $\odot$ is a transitive action of the group $G$ on $N$ verifying the equation
\begin{equation}\label{condcompatSBrcd}
g\odot (\eta\star\mu)=(g\odot \eta)\star (g\odot e_{N})^{-1}\star (g\odot \mu),
\end{equation}
for all $g\in G$ and $\eta,\mu\in N$, where $e_{N}$ denotes the unit for the group $(N,\star)$. In this initial paper, the authors emphasize in the fact that skew braces can be viewed as particular cases of skew bracoids, and prove that there exists a correspondence between Hopf-Galois structures on finite separable extensions of fields and finite skew bracoids.

Therefore, following the ideas of Angiono et al. to introduce Hopf braces from skew braces, in this paper we introduce the notion of Hopf bracoid as the quantum version of skew bracoids in a braided monoidal category ${\sf C}$. 

After an opening Section \ref{sect1} of preliminaries in which we remember the main properties of Hopf braces, Section \ref{sect2} is dedicated to introduce the category of Hopf bracoids, denoted by {\sf HBrcd}. Thereupon we present two different non-isomorphic examples of Hopf bracoids induced by a Hopf brace structure (see Examples \ref{example1} and \ref{example2}) and we also prove that generalized skew bracoids (skew bracoids for which $\odot$ is not required to be transitive) and the subcategory of all pointed and cosemisimple Hopf bracoids in ${}_{\mathbb{K}}{\sf Vect}$ such that the action is a coalgebra morphism are equivalent (see Corollary \ref{corollaryequiv}), which suppose the generalization for Hopf bracoids of this classical result in the theory of Hopf algebras. Moreover, in a Hopf brace structure it is well-known that $H_{1}$ is a left $H_{2}$-module algebra  and the action between them is also a coalgebra morphism when we work under cocommutativity assumption (see \cite[Lemma 1.8, Lemma 2.2]{AGV}). In our paper Theorem \ref{GammaB modulo algebra} and Lemma \ref{GammaBcoalg} generalize the results just mentioned to the Hopf bracoid setting. To finish this section, we prove that, when the base category is symmetric, {\sf HBrcd} is symmetric monoidal (see Theorem \ref{monoidalHBrcd}).

In \cite{AGV} it is emphasized that invertible 1-cocycles, which are no more than coalgebra isomorphisms satisfying a weaker condition than being algebra morphisms, have a hard connection with Hopf braces. To be more precise, for a fixed Hopf algebra $H_{1}$, the category of Hopf braces $\mathbb{H}=(H_{1},H_{2})$ and the category of invertible 1-cocycles $\pi\colon H_{1}\rightarrow B$ are equivalent. This categorical equivalence remains valid without the hypothesis of $H_{1}$ fixed as can be consulted in \cite[Theorem 3.2]{GONROD}. Therefore, Section \ref{sect3} of this paper contains the generalization of this result for Hopf bracoids. At first we set the notion of 1-cocycle and, after that, we construct a functor $F$ from the category of 1-cocycles to the category of Hopf bracoids (see Theorem \ref{funtor F}) and another $G$ in the opposite direction (see Theorem \ref{functorG}). This section is concluded by proving that the correspondence given by functors $F$ and $G$ induce a categorical isomorphism between certain subcategories of 1-cocycles and Hopf bracoids, as can be seen in Theorem \ref{mainiso}.

\section{Preliminaries on Hopf algebras and Hopf braces}\label{sect1}
Throughout this paper we are going to denote by ${\sf  C}$ a strict braided monoidal category with tensor product $\ot$, unit object $K$ and braiding $c$. 

As can be found in \cite{Mac}, a monoidal category is a category ${\sf  C}$ together with a functor $\ot :{\sf  C}\times {\sf  C}\rightarrow {\sf  C}$, called tensor product, an object $K$ of ${\sf C}$, called the unit object, and  families of natural isomorphisms 
$$a_{M,N,P}:(M\ot N)\ot P\rightarrow M\ot (N\ot P),\;\;\;r_{M}:M\ot K\rightarrow M, \;\;\; l_{M}:K\ot M\rightarrow M,$$
in ${\sf  C}$, called  associativity, right unit and left unit constraints, respectively, which satisfy the Pentagon Axiom and the Triangle Axiom, i.e.,
$$a_{M,N, P\ot Q}\co a_{M\ot N,P,Q}= (id_{M}\ot a_{N,P,Q})\co a_{M,N\ot P,Q}\co (a_{M,N,P}\ot id_{Q}),$$
$$(id_{M}\ot l_{N})\co a_{M,K,N}=r_{M}\ot id_{N},$$
where for each object $X$ in ${\sf  C}$, $id_{X}$ denotes the identity morphism of $X$. A monoidal category is called strict if the previous constraints are identities. It is an important result (see for example \cite{K}) that every non-strict monoidal category is monoidal equivalent to a strict one, so the strict character can be assumed without loss of generality. Then, results proved in a strict setting hold for every non-strict monoidal category that include, between others, the category ${}_{\mathbb{K}}\sf{Vect}$ of vector spaces over a field ${\mathbb K}$ and the category ${}_{R}\sf{Mod}$ of left modules over a commutative ring $R$.
For
simplicity of notation, given objects $M$, $N$, $P$ in ${\sf  C}$ and a morphism $f:M\rightarrow N$, in most cases we will write $P\ot f$ for $id_{P}\ot f$ and $f \ot P$ for $f\ot id_{P}$.

A braiding for a strict monoidal category ${\sf  C}$ is a natural family of isomorphisms 
$$c_{M,N}:M\ot N\rightarrow N\ot M$$ subject to the conditions 
$$
c_{M,N\ot P}= (N\ot c_{M,P})\co (c_{M,N}\ot P),\;\;
c_{M\ot N, P}= (c_{M,P}\ot N)\co (M\ot c_{N,P}).
$$

A strict braided monoidal category ${\sf  C}$ is a strict monoidal category with a braiding. These categories were introduced by Joyal and Street (see \cite{JS2}) motivated by the theory of braids and links in topology. Note that, as a consequence of the definition, the equalities $c_{M,K}=c_{K,M}=id_{M}$ hold, for all object  $M$ of ${\sf  C}$. 

If the braiding satisfies that  $c_{N,M}\co c_{M,N}=id_{M\ot N},$ for all $M$, $N$ in ${\sf  C}$, we will say that ${\sf C}$  is symmetric. In this case, we call the braiding $c$ a symmetry for the category ${\sf  C}$.

In what follows we will remember some basic definitions that will be useful in the paper.

\begin{definition}
{\rm 
An algebra in ${\sf  C}$ is a triple $A=(A, \eta_{A},
\mu_{A})$ where $A$ is an object in ${\sf  C}$ and
 \mbox{$\eta_{A}:K\rightarrow A$} (unit), $\mu_{A}:A\otimes A
\rightarrow A$ (product) are morphisms in ${\sf  C}$ such that
$\mu_{A}\circ (A\otimes \eta_{A})=id_{A}=\mu_{A}\circ
(\eta_{A}\otimes A)$, $\mu_{A}\circ (A\otimes \mu_{A})=\mu_{A}\circ
(\mu_{A}\otimes A)$. Given two algebras $A= (A, \eta_{A}, \mu_{A})$
and $B=(B, \eta_{B}, \mu_{B})$, a morphism  $f:A\rightarrow B$ in {\sf  C} is an algebra morphism if $\mu_{B}\circ (f\otimes f)=f\circ \mu_{A}$, $ f\circ
\eta_{A}= \eta_{B}$. 

If  $A$, $B$ are algebras in ${\sf  C}$, the tensor product
$A\otimes B$ is also an algebra in
${\sf  C}$ where
$\eta_{A\otimes B}=\eta_{A}\otimes \eta_{B}$ and $\mu_{A\otimes
	B}=(\mu_{A}\otimes \mu_{B})\circ (A\otimes c_{B,A}\otimes B).$
}
\end{definition}

\begin{definition}
{\rm 
A coalgebra  in ${\sf  C}$ is a triple ${D} = (D,
\varepsilon_{D}, \delta_{D})$ where $D$ is an object in ${\sf
C}$ and $\varepsilon_{D}: D\rightarrow K$ (counit),
$\delta_{D}:D\rightarrow D\otimes D$ (coproduct) are morphisms in
${\sf  C}$ such that $(\varepsilon_{D}\otimes D)\circ
\delta_{D}= id_{D}=(D\otimes \varepsilon_{D})\circ \delta_{D}$,
$(\delta_{D}\otimes D)\circ \delta_{D}=
 (D\otimes \delta_{D})\circ \delta_{D}.$ If ${D} = (D, \varepsilon_{D},
 \delta_{D})$ and
${ E} = (E, \varepsilon_{E}, \delta_{E})$ are coalgebras, a morphism 
$f:D\rightarrow E$ in  {\sf  C} is a coalgebra morphism if $(f\otimes f)\circ
\delta_{D} =\delta_{E}\circ f$, $\varepsilon_{E}\circ f
=\varepsilon_{D}.$ 

Given  $D$, $E$ coalgebras in ${\sf  C}$, the tensor product $D\otimes E$ is a
coalgebra in ${\sf  C}$ where $\varepsilon_{D\otimes
E}=\varepsilon_{D}\otimes \varepsilon_{E}$ and $\delta_{D\otimes
E}=(D\otimes c_{D,E}\otimes E)\circ( \delta_{D}\otimes \delta_{E}).$
}
\end{definition}

\begin{definition}
Consider $(C,\varepsilon_{C},\delta_{C})$ a coalgebra in ${}_{\mathbb{K}}\sf{Vect}$. An element $c\in C$ is a group-like element in $C$ if $\delta_{C}(c)=c\otimes c$ and $\varepsilon_{C}(c)=1_{\mathbb{K}}$. 

In what follows, ${\sf G}(C)$ will denote the set of group-like elements in $C$, which is a subcoalgebra of the coalgebra $C$.
\end{definition}

\begin{example}\label{freevectspacecoalgebra}
Given a set $S$, the free $\mathbb{K}$-vector space $\mathbb{K}[S]\coloneqq\bigoplus_{s\in S}\mathbb{K}s$ has a natural coalgebra structure in ${}_{\mathbb{K}}\sf{Vect}$ whose coproduct is defined by $\delta_{\mathbb{K}[S]}(s)=s\otimes s$ and $\varepsilon_{\mathbb{K}[S]}(s)=1_{\mathbb{K}}$ for all $s\in S$. For this coalgebra, it is straightforward to compute that ${\sf G}(\mathbb{K}[S])=S$. In the particular case that the set is a group $(G,\cdot)$, ${\sf G}(\mathbb{K}[G])=G$ which is also a group.
\end{example}

\begin{definition}
A coalgebra $(C,\varepsilon_{C},\delta_{C})$ in ${}_{\mathbb{K}}\sf{Vect}$ is pointed if its simple subcoalgebras are one-dimensional. As a characterization of this property, $C$ is pointed if and only if its coradical $C_{0}$ (the sum of all simple subcoalgebras) is $\mathbb{K}[{\sf G}(C)]$ (see \cite[Section 3.4]{RADHA}).
\end{definition}
\begin{definition}
A coalgebra $(C,\varepsilon_{C},\delta_{C})$ in ${}_{\mathbb{K}}\sf{Vect}$ is cosemisimple if $C_{0}=C$.
\end{definition}

So, putting together the characterization of pointed coalgebra with cosemisimplicity, we obtain that a coalgebra $(C,\varepsilon_{C},\delta_{C})$ is pointed and cosemisimple if and only if 
\begin{equation}\label{pointed+cosemis}
C=\mathbb{K}[{\sf G}(C)].
\end{equation}

\begin{definition}
	{\rm 
 Let ${D} = (D, \varepsilon_{D},
\delta_{D})$ be a coalgebra and $A=(A, \eta_{A}, \mu_{A})$ an
algebra in $\sf{C}$. By ${\mathcal  H}(D,A)$ we denote the set of morphisms
$f:D\rightarrow A$ in ${\sf  C}$. With the convolution operation
$f\ast g= \mu_{A}\circ (f\otimes g)\circ \delta_{D}$, ${\mathcal  H}(D,A)$ is an algebra where the unit element is $\eta_{A}\circ \varepsilon_{D}=\varepsilon_{D}\otimes \eta_{A}$.
}
\end{definition}

\begin{definition}
{\rm 
 Let  $A$ be an algebra. The pair
$(M,\varphi_{M})$ is a left $A$-module if $M$ is an object in
${\sf  C}$ and $\varphi_{M}:A\otimes M\rightarrow M$ is a morphism
in ${\sf  C}$ satisfying $\varphi_{M}\circ(
\eta_{A}\ot M)=id_{M}$, $\varphi_{M}\circ (A\ot \varphi_{M})=\varphi_{M}\circ
(\mu_{A}\ot M)$. Given two left ${A}$-modules $(M,\varphi_{M})$
and $(N,\varphi_{N})$, $f:M\rightarrow N$ is a morphism of left
${A}$-modules if $\varphi_{N}\circ (A\ot f)=f\circ \varphi_{M}$.  

Due to being the  composition of morphisms of left $A$-modules a morphism of left $A$-modules, left $A$-modules form a category that we will denote by $\;_{\sf A}${\sf Mod}.

}

\end{definition}

\begin{definition}
{\rm 
We say that $X$ is a
bialgebra  in ${\sf  C}$ if $(X, \eta_{X}, \mu_{X})$ is an
algebra, $(X, \varepsilon_{X}, \delta_{X})$ is a coalgebra, and
$\varepsilon_{X}$ and $\delta_{X}$ are algebra morphisms
(equivalently, $\eta_{X}$ and $\mu_{X}$ are coalgebra morphisms). Moreover, if there exists a morphism $\lambda_{X}:X\rightarrow X$ in ${\sf  C}$,
called the antipode of $X$, satisfying that $\lambda_{X}$ is the inverse of $id_{X}$ in ${\mathcal  H}(X,X)$, i.e., 
\begin{equation}
\label{antipode}
id_{X}\ast \lambda_{X}= \eta_{X}\circ \varepsilon_{X}= \lambda_{X}\ast id_{X},
\end{equation}
we say that $X$ is a Hopf algebra. A morphism of Hopf algebras is an algebra-coalgebra morphism. Note that, if $f:X\rightarrow Y$ is a Hopf algebra morphism the following equality holds:
\begin{equation}
\label{morant}
\lambda_{Y}\co f=f\co \lambda_{X}.
\end{equation} 

With the composition of morphisms in {\sf C} we can define a category whose objects are Hopf algebras  and whose morphisms are morphisms of Hopf algebras. We denote this category by ${\sf  Hopf}$.

\begin{example}\label{functor L y G}
Note that $X$ is a Hopf algebra in $\sf{Set}$ if and only if $X$ is a group. So, if $\mathbb{K}$ is a field and $(G,\cdot)$ is a group with unit $e_{G}$, the free $\mathbb{K}$-vector space $\mathbb{K}[G]$ is a Hopf algebra in ${}_{\mathbb{K}}\sf{Vect}$, where the coalgebra structure is the one described in Example \ref{freevectspacecoalgebra}, and the algebra structure is given by 
\[
\eta_{\mathbb{K}[G]}(1_{\mathbb{K}})=e_{G},\;\;\mu_{\mathbb{K}[G]}(a\otimes b)=a\cdot b,
\]
with antipode $\lambda_{\mathbb{K}[G]}(a)=a^{-1}$, for all $a,b\in G$. Moreover, if $f\colon (G,\cdot)\rightarrow (N,\ast)$ is a group morphism, then $\mathbb{K}[f]\colon \mathbb{K}[G]\rightarrow \mathbb{K}[N]$, which is the liner extension of $f$ to $\mathbb{K}[G]$, is a Hopf algebra morphism in ${}_{\mathbb{K}}\sf{Vect}$. Therefore, if $\sf{Grp}$ denotes the category of groups and ${}_{\mathbb{K}}\sf{Hopf}$, the category of Hopf algebras in ${}_{\mathbb{K}}\sf{Vect}$, then there exists a functor
\[\sf{L}\colon \sf{Grp}\longrightarrow {}_{\mathbb{K}}\sf{Hopf}\]
acting on objects by ${\sf L}((G,\cdot))=\mathbb{K}[G]$ and on morphisms by ${\sf L}(f)=\mathbb{K}[f]$.

Conversely, given $X=(X,\eta_{X},\mu_{X},\varepsilon_{X},\delta_{X},\lambda_{X})$ a Hopf algebra in ${}_{\mathbb{K}}\sf{Vect}$, ${\sf G}(X)$ is a group (see \cite[Proposition III.3.7]{K}) whose product, unit element and inverse are defined as follows:
\[
a \diamond b\coloneqq \mu_{X}(a\otimes b),\;\; e_{{\sf G}(X)}\coloneqq \eta_{X}(1_{\mathbb{K}}),\;\; a^{-1}\coloneqq \lambda_{X}(a),
\]
for all $a\in {\sf G}(X)$. In addition, if $f\colon X\rightarrow H$ is a Hopf algebra morphism in ${}_{\mathbb{K}}\sf{Vect}$, then the restriction of $f$ to ${\sf G}(X)$ is a group morphism whose image lies in ${\sf G}(H)$. In conclusion, there exists a functor 
\[{\sf G}\colon {}_{\mathbb{K}}\sf{Hopf}\longrightarrow \sf{Grp}\]
which acts on objects sending every Hopf algebra $X$ in ${}_{\mathbb{K}}\sf{Vect}$ to its set of group-like elements, ${\sf G}(X)$, and on morphisms sending every Hopf algebra morphism $f\colon X\rightarrow H$ in ${}_{\mathbb{K}}\sf{Vect}$ to its restriction to the group-like elements.

The previous correspondence gives rise to an adjoint pair, ${\sf L}\dashv{\sf G}$, which induce a categorical equivalence between ${\sf Grp}$ and the full subcategory of ${}_{\mathbb{K}}\sf{Hopf}$ whose objects are pointed and cosemisimple Hopf algebras in ${}_{\mathbb{K}}\sf{Vect}$.
\end{example}

A Hopf algebra is commutative if $\mu_{X}\co c_{X,X}=\mu_{X}$ and cocommutative if $c_{X,X}\co \delta_{X}=\delta_{X}.$ It is easy to see that in both cases $\lambda_{X}\circ \lambda_{X} =id_{X}$.
}
\end{definition}

If $X$ is a Hopf algebra, significant properties of its antipode, $\lambda_{X}$, are the following:  It is antimultiplicative and anticomultiplicative 
\begin{equation}
\label{a-antip}
\lambda_{X}\co \mu_{X}=  \mu_{X}\co (\lambda_{X}\ot \lambda_{X})\co c_{X,X},\;\;\;\; \delta_{X}\co \lambda_{X}=c_{X,X}\co (\lambda_{X}\ot \lambda_{X})\co \delta_{X}, 
\end{equation}
and leaves the unit and counit invariant, i.e., 
\begin{equation}
\label{u-antip}
\lambda_{X}\co \eta_{X}=  \eta_{X},\;\; \varepsilon_{X}\co \lambda_{X}=\varepsilon_{X}.
\end{equation}
So, it is a direct consequence of these identities that, if $X$ is commutative, then $\lambda_{X}$ is an algebra morphism and, if $X$ is cocommutative, then $\lambda_{X}$ is a coalgebra morphism.

In the following definition we recall the notion of left module (co)algebra.

\begin{definition}
{\rm 
Let $X$ be a Hopf algebra. An algebra $A$  is said to be a left $X$-module algebra if $(A, \varphi_{A})$ is a left $X$-module and $\eta_{A}$, $\mu_{A}$ are morphisms of left $X$-modules, i.e.,
\begin{equation}
\label{mod-alg1}
\varphi_{A}\circ (X\otimes \eta_{A})=\varepsilon_{X}\otimes \eta_{A},
\end{equation}
\begin{equation}\label{mod-alg2}
\varphi_{A}\circ (X\otimes \mu_{A})=\mu_{A}\circ \varphi_{A\otimes A},
\end{equation}
where  $\varphi_{A\otimes A}=(\varphi_{A}\otimes \varphi_{A})\circ (X\otimes c_{X,A}\otimes A)\circ (\delta_{X}\otimes A\otimes A)$ is the left action on $A\otimes A$. 
}
\end{definition}
\begin{definition}
Let $X$ be a Hopf algebra. A coalgebra $D$ is said to be a left $X$-module coalgebra if $(D,\varphi_{D})$ is a left $X$-module and $\varepsilon_{D}$, $\delta_{D}$ are morphisms of left $X$-modules, in other words, the following equalities hold:
\begin{equation}\label{mod-coalg1}
\varepsilon_{D}\circ\varphi_{D}=\varepsilon_{X}\otimes\varepsilon_{D},
\end{equation}
\begin{equation}\label{mod-coalg2}
\delta_{D}\circ\varphi_{D}=\varphi_{D\otimes D}\circ(X\otimes\delta_{D}).
\end{equation}
Equivalently, $(D,\varphi_{D})$ is a left $X$-module coalgebra if and only if $\varphi_{D}$ is a coalgebra morphism.
\end{definition}

In \cite{AGV}, Angiono, Galindo and Vendramin introduced a generalisation of the Hopf algebra structure, the so called Hopf braces. In the braided setting the definition of Hopf brace is the following:
\begin{definition}
\label{H-brace}
{\rm Let $H=(H, \varepsilon_{H}, \delta_{H})$ be a coalgebra in {\sf C}. Let's assume that there are two algebra structures $(H, \eta_{H}^1, \mu_{H}^1)$, $(H, \eta_{H}^2, \mu_{H}^2)$ defined on $H$ and suppose that there exist two endomorphism of $H$ denoted by $\lambda_{H}^{1}$ and $\lambda_{H}^{2}$. We will say that 
$$(H, \eta_{H}^{1}, \mu_{H}^{1}, \eta_{H}^{2}, \mu_{H}^{2}, \varepsilon_{H}, \delta_{H}, \lambda_{H}^{1}, \lambda_{H}^{2})$$
is a Hopf brace in {\sf C} if:
\begin{itemize}
\item[(i)] $H_{1}=(H, \eta_{H}^{1}, \mu_{H}^{1},  \varepsilon_{H}, \delta_{H}, \lambda_{H}^{1})$ is a Hopf algebra in {\sf C}.
\item[(ii)] $H_{2}=(H, \eta_{H}^{2}, \mu_{H}^{2},  \varepsilon_{H}, \delta_{H}, \lambda_{H}^{2})$ is a Hopf algebra in {\sf C}.
\item[(iii)] The  following equality holds:
$$\mu_{H}^{2}\co (H\ot \mu_{H}^{1})=\mu_{H}^{1}\co (\mu_{H}^{2}\ot \Gamma_{H_{1}} )\co (H\ot c_{H,H}\ot H)\co (\delta_{H}\ot H\ot H),$$
\end{itemize}
where  $$\Gamma_{H_{1}}\coloneqq\mu_{H}^{1}\co (\lambda_{H}^{1}\ot \mu_{H}^{2})\co (\delta_{H}\ot H).$$

Following \cite{RGON}, a Hopf brace will be denoted by ${\mathbb H}=(H_{1}, H_{2})$ or in a simpler way by $\mathbb{H}$.

}
\end{definition}
\begin{definition}
%\label{H-coco}
 If  ${\mathbb H}$ is a Hopf brace in {\sf C}, we will say that ${\mathbb H}$ is cocommutative if $\delta_{H}=c_{H,H}\circ \delta_{H}$, i.e., if $H_{1}$ and $H_{2}$ are cocommutative Hopf algebras in {\sf C}.
 
%Note that by \cite[Corollary 5]{Sch}, if $H$ is a  cocommutative Hopf algebra  in the  braided monoidal category {\sf C}, the identity 
%\begin{equation}
%\label{ccb}
%c_{H,H}\circ c_{H,H}=id_{H\otimes H} 
%\end{equation}
%holds.
\end{definition}

\begin{definition}
%\label{mor}
{\rm  Given two Hopf braces ${\mathbb H}$  and  ${\mathbb B}$ in {\sf C}, a morphism $f$ in {\sf C} between the two underlying objects is called a morphism of Hopf braces if both $f:H_{1}\rightarrow B_{1}$ and $f:H_{2}\rightarrow B_{2}$ are Hopf algebra morphisms.
		
Hopf braces together with morphisms of Hopf braces form a category which we denote by {\sf HBr}. Moreover, cocommutative Hopf braces constitute a full subcategory of $\sf{HBr}$ which we will denote by $\sf{cocHBr}$.
		
}
\end{definition}

Between the properties of these objects the following are noteworthy: Let ${\mathbb H}$  be a Hopf brace in {\sf C}, then 
\begin{equation}
%\label{eb1}
\eta_{H}^{1}=\eta_{H}^2
 \end{equation}
holds and, by \cite[Lemma 1.7]{AGV},  in this braided setting  the equality
\begin{equation}
\label{agv1}
\Gamma_{H_{1}}\circ (H\otimes \lambda_{H}^1)=\mu_{H}^{1}\circ ((\lambda_{H}^1\circ \mu_{H}^{2})\otimes H)\circ (H\otimes c_{H,H}) \circ (\delta_{H}\otimes H)
\end{equation}
also holds. Moreover, in our braided context \cite[Lemma 1.8]{AGV} and \cite[Remark 1.9]{AGV} hold and then we have that the algebra $(H,\eta_{H}^{1}, \mu_{H}^{1})$ is a left $H_{2}$-module algebra with action $\Gamma_{H_{1}}$ and $\mu_{H}^2$ admits the following expression:
\begin{equation}
\label{eb2}
\mu_{H}^2=\mu_{H}^{1}\circ (H\otimes \Gamma_{H_{1}})\circ (\delta_{H}\otimes H). 
\end{equation}
In addition, by \cite[Lemma 2.2]{AGV}, $\Gamma_{H_{1}}$ is a coalgebra morphism when $\mathbb{H}$ is cocommutative.

\section{Hopf bracoids}\label{sect2}

In this section, we are going to introduce the notion of Hopf bracoid in a braided monoidal setting. Hopf bracoids are the quantum version of skew bracoids that have been introduced by Martin-Lyons and Paul J. Truman in \cite{MLT}. Most of the properties that we can check on \cite[Section 1]{MLT} have an analogous in this new context.

We will define a Hopf bracoid in $\sf{C}$ as follows.
\begin{definition}\label{HBrcd} Let $H= (H,\eta_{H},\mu_{H},\varepsilon_{H},\delta_{H},\lambda_{H})$ and $B=(B,\eta_{B},\mu_{B},\varepsilon_{B},\delta_{B},\lambda_{B})$ be Hopf algebras in $\sf{C}$. We will say that a triple $(H,B,\varphi_{B})$ is a Hopf bracoid in $\sf{C}$ if $(B,\varphi_{B})$ is a left $H$-module such that the following compatibility condition holds:
\begin{equation}\label{Eq. Bracoid compatibility}
\varphi_{B}\circ(H\otimes\mu_{B})=\mu_{B}\circ (\varphi_{B}\otimes\Phi_{B})\circ(H\otimes c_{H,B}\otimes B)\circ(\delta_{H}\otimes B\otimes B),
\end{equation}
where $\Phi_{B}\coloneqq \mu_{B}\circ((\lambda_{B}\circ u_{\varphi_{B}})\otimes\varphi_{B})\circ(\delta_{H}\otimes B)$ and $u_{\varphi_{B}}\coloneqq\varphi_{B}\circ(H\otimes\eta_{B})$.
\end{definition}

\begin{remark}
Note that Hopf bracoids are a generalization of module algebras in a certain sense. That is because if $u_{\varphi_{B}}=\varepsilon_{H}\otimes\eta_{B}$, then \eqref{Eq. Bracoid compatibility} implies that $\mu_{B}$ is a morphism of left $H$-modules and, as a consequence, $(B,\varphi_{B})$ is a left $H$-module algebra.
\end{remark}
\begin{definition}
Let $(H,B,\varphi_{B})$ and $(H',B',\varphi_{B'})$ be Hopf bracoids in $\sf{C}$. A morphism of Hopf bracoids in $\sf{C}$ between $(H,B,\varphi_{B})$ and $(H',B',\varphi_{B'})$ is a pair $(f,g)$ where $f\colon H\rightarrow H'$ and $g\colon B\rightarrow B'$ are morphisms of Hopf algebras in $\sf{C}$ such that the following condition holds:
\begin{equation}\label{Cond. Mor bracoid}
g\circ\varphi_{B}=\varphi_{B'}\circ(f\otimes g).
\end{equation}

Therefore, Hopf bracoids in $\sf{C}$ give rise to a category whose morphisms are the ones defined in the previous definition. We will denote it by $\hbr$. 
\end{definition}
\begin{example}\label{example1}
If $\mathbb{H}=(H_{1},H_{2})$ is a Hopf brace in $\sf{C}$, then  $(H_{2},H_{1},\mu_{H}^{2})$ is a Hopf bracoid in $\sf{C}$. So, there exists a functor ${\sf T}\colon\textnormal{\sf{HBr}}\longrightarrow \hbr$ acting on objects by ${\sf T}(\mathbb{H})=(H_{2},H_{1},\mu_{H}^{2})$ and on morphisms by ${\sf T}(f)=(f,f)$. We will call this functor the trivial functor from Hopf braces to Hopf bracoids.
\end{example}
\begin{example}\label{example2}
Consider $\mathbb{H}=(H_{1},H_{2})$ a Hopf brace in ${\sf C}$. Thanks to general properties of Hopf braces, $(H_{1},\Gamma_{H_{1}})$ is a left $H_{2}$-module algebra. Then, the triple $(H_{2},H_{1},\Gamma_{H_{1}})$ is a Hopf bracoid in ${\sf C}$. Indeed, we have that
\begin{align*}
u_{\Gamma_{H_{1}}}=\Gamma_{H_{1}}\circ(H\otimes\eta_{H})=\varepsilon_{H}\otimes\eta_{H}\;\footnotesize\textnormal{(by the condition of morphism of left $H_{2}$-modules for $\eta_{H}$)}
\end{align*}
and then, thanks to the above equality, \eqref{u-antip} and (co)unit properties,
\begin{align}\label{PhiH1=GammaH1}
\Phi_{H_{1}}=\Gamma_{H_{1}}.
\end{align}
This implies that
\begin{align*}
&\mu_{H}^{1}\circ(\Gamma_{H_{1}}\otimes\Phi_{H_{1}})\circ(H\otimes c_{H,H}\otimes H)\circ(\delta_{H}\otimes H\otimes H)\\=&\mu_{H}^{1}\circ(\Gamma_{H_{1}}\otimes\Gamma_{H_{1}})\circ(H\otimes c_{H,H}\otimes H)\circ(\delta_{H}\otimes H\otimes H)\;\footnotesize\textnormal{(by \eqref{PhiH1=GammaH1})}\\=&\Gamma_{H_{1}}\circ(H\otimes\mu_{H}^{1})\;\footnotesize\textnormal{(by the condition of morphism of left $H_{2}$-modules for $\mu_{H}^{1}$)}.
\end{align*}
In addition, if $\mathbb{H}$ and $\mathbb{B}$ are Hopf braces and $f$ is a morphism between them, then $(f,f)$ is a morphism in ${\sf HBrcd}$ between $(H_{2},H_{1},\Gamma_{H_{1}})$ and $(B_{2},B_{1},\Gamma_{B_{1}})$ because \eqref{Cond. Mor bracoid} holds. Indeed,
\begin{align*}
&f\circ\Gamma_{H_{1}}\\=&\mu_{B}^{1}\circ((f\circ\lambda_{H}^{1})\otimes (f\circ\mu_{H}^{2}))\circ(\delta_{H}\otimes H)\;\footnotesize\textnormal{(by the condition of algebra morphism for $f\colon H_{1}\rightarrow B_{1}$)}\\=&\mu_{B}^{1}\circ((\lambda_{B}^{1}\circ f)\otimes(\mu_{B}^{2}\circ(f\otimes f)))\circ(\delta_{H}\otimes H)\;\footnotesize\textnormal{(by the condition of algebra morphism for $f\colon H_{2}\rightarrow B_{2}$ and \eqref{morant})}\\=&\mu_{B}^{1}\circ(\lambda_{B}^{1}\otimes\mu_{B}^{2})\circ(\delta_{B}\otimes H)\circ(f\otimes f)\;\footnotesize\textnormal{(by the condition of coalgebra morphism for $f$)}\\=&\Gamma_{B_{1}}\circ(f\otimes f).
\end{align*}
As a conclusion, there exists a functor $T'\colon {\sf HBr}\longrightarrow {\sf HBrcd}$ acting on objects by $T'(\mathbb{H})=(H_{2},H_{1},\Gamma_{H_{1}})$ and on morphisms by $T'(f)=(f,f)$.
\end{example}
\begin{example}\label{example functor P}
Following \cite{MLT}, a (left) skew bracoid is a 5-tuple $(G,\cdot,N,\ast,\odot)$ such that $(G,\cdot)$ and $(N,\ast)$ are groups, $\odot\colon G\times N\rightarrow N$ is a transitive action of $(G,\cdot)$ on $N$ and the following compatibility condition holds:
\begin{equation}\label{skbracoidcompatibility}
g\odot (n\ast n')=(g\odot n)\ast(g\odot e_{N})^{-1}\ast(g\odot n')
\end{equation}
for all $g\in G$ and $n,n'\in N$, and where $e_{N}$ denotes the unit element for $(N,\ast)$. 

If in the definition of  skew bracoid we do not assume that the action $\odot$ is transitive we have the definition of generalized skew bracoid. Given $(G,\cdot,N,\ast,\odot)$ and $(G',\cdot',N',\ast',\odot')$ generalized skew bracoids, a pair $(f,g)$, where $f\colon (G,\cdot)\rightarrow (G',\cdot^{\prime})$ and $g\colon (N,\ast)\rightarrow (N',\ast')$ are group morphisms, is said to be a morphism of generalized skew bracoids if the identity 
\begin{equation}\label{morgenskbracoid}
g\circ\odot=\odot'\circ(f\otimes g)
\end{equation}
holds. So, generalized skew bracoids constitute a category that we will denote by $\textnormal{\texttt{g}}\sf{SkBrcd}$ for which skew bracoids are a full subcategory denoted by $\sf{SkBrcd}$. As a conclusion,  $(G,\cdot,N,\ast,\odot)$ is a generalized skew bracoid if and only if it is a Hopf bracoid in {\sf Set}.

Let $\mathbb{K}$ be a field. If $(G,\cdot,N,\ast,\odot)$ is an object in $\textnormal{\texttt{g}}\sf{SkBrcd}$, then $(\mathbb{K}[G],\mathbb{K}[N],\varphi_{\odot})$ is a Hopf bracoid in ${}_{\mathbb{K}}\sf{Vect}$, where $\mathbb{K}[G]={\sf L}((G,\cdot))$ and $\mathbb{K}[N]={\sf L}((N,\ast))$, being $\sf{L}$ the functor defined in Example \ref{functor L y G}, and $\varphi_{\odot}$ is the linear extension of $\odot$ to $\mathbb{K}[G]\otimes\mathbb{K}[N]$. Also if $(f,g)\colon (G,\cdot,N,\ast,\odot)\rightarrow (G',\cdot',N',\ast',\odot')$ is a morphism in $\textnormal{\texttt{g}}\sf{SkBrcd}$, then $(\mathbb{K}[f],\mathbb{K}[g])$ is a morphism of Hopf bracoids in ${}_{\mathbb{K}}\sf{Vect}$ between $(\mathbb{K}[G],\mathbb{K}[N],\varphi_{\odot})$ and $(\mathbb{K}[G'],\mathbb{K}[N'],\varphi_{\odot'})$. Therefore, denoting by ${}_{\mathbb{K}}\sf{HBrcd}$ the category of Hopf bracoids in ${}_{\mathbb{K}}\sf{Vect}$, there exists a functor 
\[\sf{P}\colon \textnormal{\texttt{g}}\sf{SkBrcd}\longrightarrow {}_{\mathbb{K}}\sf{HBrcd}\]
acting on objects by ${\sf P}((G,\cdot,N,\ast,\odot))=({\sf L}((G,\cdot)),{\sf L}((N,\ast)),\varphi_{\odot})$ and on morphisms by ${\sf P}((f,g))=(\mathbb{K}[f],\mathbb{K}[g])$.
\end{example}

\begin{example}\label{example functor R}
Let $\mathbb{K}$ be a field and $(H,B,\varphi_{B})$ a Hopf bracoid in ${}_{\mathbb{K}}\sf{Vect}$ such that $\varphi_{B}$ is a coalgebra morphism. Let's prove that $({\sf G}(H),{\sf G}(B),\circledast)$ is a generalized skew bracoid, where ${\sf G}$ is the group-like functor defined in Example \ref{functor L y G} and $\circledast$ is an action of ${\sf G}(H)$ over ${\sf G}(B)$ defined by $h\circledast b\coloneqq \varphi_{B}(h\otimes b)$ for all $h\in {\sf G}(H)$ and $b\in {\sf G}(B)$. First of all we have to see that $h\circledast b$ is a group-like element in $B$ for all $h\in {\sf G}(H)$ and $b\in {\sf G}(B)$. Indeed,
\begin{align*}
&\delta_{B}(h\circledast b)=\delta_{B}(\varphi_{B}(h\otimes b))\\=&((\varphi_{B}\otimes\varphi_{B})\circ(H\otimes c_{H,B}\otimes B)\circ(\delta_{H}\otimes \delta_{B}))(h\otimes b)\;\footnotesize\textnormal{(by the condition of coalgebra morphism for $\varphi_{B}$)}\\=&\varphi_{B}(h\otimes b)\otimes\varphi_{B}(h\otimes b)\;\footnotesize\textnormal{(by $h\in{\sf G}(H)$ and $b\in{\sf G}(B)$)}\\=&(h\circledast b)\otimes(h\circledast b)
\end{align*}
and, moreover,
\begin{align*}
&\varepsilon_{B}(h\circledast b)=\varepsilon_{B}(\varphi_{B}(h\otimes b))\\=&\varepsilon_{H}(h)\otimes\varepsilon_{B}(b)\;\footnotesize\textnormal{(by the condition of coalgebra morphism for $\varphi_{B}$)}\\=&1_{\mathbb{K}}\otimes 1_{\mathbb{K}}\;\footnotesize\textnormal{(by $h\in{\sf G}(H)$ and $b\in{\sf G}(B)$)}\\=&1_{\mathbb{K}}.
\end{align*}
So, $h\circledast b\in{\sf G}(B)$ for all $h\in{\sf G}(H)$ and $b\in{\sf G}(B)$. We already know that ${\sf G}(H)$ and ${\sf G}(B)$ are groups and, in addition, $\circledast$ is an action because $\varphi_{B}$ is also an action, so to conclude that $({\sf G}(H),{\sf G}(B),\circledast)$ is a generalized skew bracoid it is enough to compute that \eqref{skbracoidcompatibility} holds, what follows directly thanks to \eqref{Eq. Bracoid compatibility}. Moreover, due to the functoriality of ${\sf G}$, given $(f,g)\colon (H,B,\varphi_{B})\rightarrow (H',B',\varphi_{B'})$ a morphism in ${}_{\mathbb{K}}{\sf HBrcd}$ between Hopf bracoids satisfying that $\varphi_{B}$ and $\varphi_{B^{\prime}}$ are coalgebra morphisms, $({\sf G}(f),{\sf G}(g))$ is a morphism in $\textnormal{\texttt{g}}{\sf SkBrcd}$ because condition \eqref{morgenskbracoid} follows from \eqref{Cond. Mor bracoid}, where ${\sf G}(f)$ and ${\sf G}(g)$ are the restrictions of $f$ and $g$, respectively, to the group-like elements. As a conclusion, if we denote by ${}_{\mathbb{K}}\overline{{\sf HBrcd}}$ to the full subcategory of ${}_{\mathbb{K}}{\sf HBrcd}$ whose objects are Hopf bracoids $(H,B,\varphi_{B})$ in ${}_{\mathbb{K}}\sf{Vect}$ such that $\varphi_{B}$ is a coalgebra morphism, there exists a functor \[{\sf R}\colon {}_{\mathbb{K}}\overline{{\sf HBrcd}}\longrightarrow \textnormal{\texttt{g}}{\sf SkBrcd}\] acting on objects by ${\sf R}((H,B,\varphi_{B}))=({\sf G}(H),{\sf G}(B),\circledast)$ and on morphisms by ${\sf R}((f,g))=({\sf G}(f),{\sf G}(g))$.
\end{example}
\begin{remark}
Note that, if ${\sf P}$ is the functor defined in Example \ref{example functor P} and $(G,\cdot,N,\ast,\odot)$ is an object in $\textnormal{\texttt{g}}{\sf SkBrcd}$, then ${\sf P}((G,\cdot,N,\ast,\odot))=({\sf L}((G,\cdot)),{\sf L}((N,\ast)),\varphi_{\odot})$ is an object in ${}_{\mathbb{K}}\overline{{\sf HBrcd}}$, that is to say, $\varphi_{\odot}$ is a coalgebra morphism. Indeed, consider $g\in G$ and $n\in N$, we have that
\begin{align*}
&(\delta_{\mathbb{K}[N]}\circ\varphi_{\odot})(g\otimes n)=\delta_{\mathbb{K}[N]}(g\odot n)=(g\odot n)\otimes(g\odot n)\\=&((\varphi_{\odot}\otimes\varphi_{\odot})\circ(\mathbb{K}[G]\otimes c_{\mathbb{K}[G],\mathbb{K}[N]}\otimes\mathbb{K}[N])\circ(\delta_{\mathbb{K}[G]}\otimes \delta_{\mathbb{K}[N]}))(g\otimes n),
\end{align*}
and also
\begin{align*}
&(\varepsilon_{\mathbb{K}[N]}\circ\varphi_{\odot})(g\otimes n)=\varepsilon_{\mathbb{K}[N]}(g\odot n)=1_{\mathbb{K}}=(\varepsilon_{\mathbb{K}[G]}\otimes\varepsilon_{\mathbb{K}[N]})(g\otimes n).
\end{align*}
Extending previous equalities by linearity to $\mathbb{K}[G]\otimes\mathbb{K}[N]$, we obtain that $\varphi_{\odot}$ is a coalgebra morphism. So, ${\sf P}$ is a functor from $\textnormal{\texttt{g}}{\sf SkBrcd}$ to ${}_{\mathbb{K}}\overline{{\sf HBrcd}}$.
\end{remark}
\begin{theorem}\label{adjunction}
The functor ${\sf P}$ defined in Example \ref{example functor P} is left adjoint to the functor ${\sf R}$ defined in Example \ref{example functor R}.
\end{theorem}
\begin{proof}
Let $(G,\cdot,N,\ast,\odot)$ be an object in $\textnormal{\texttt{g}}{\sf SkBrcd}$ and let $(H,B,\varphi_{B})$ be an object in ${}_{\mathbb{K}}\overline{{\sf HBrcd}}$. We have to define a bijection 
\[{}^{(G,\cdot,N,\ast,\odot)}\Gamma_{(H,B,\varphi_{B})}\colon \begin{array}{c}\operatorname{Hom}_{{}_{\mathbb{K}}\overline{{\sf HBrcd}}}({\sf P}((G,\cdot,N,\ast,\odot)),(H,B,\varphi_{B}))\\\downarrow\\ \operatorname{Hom}_{\textnormal{\texttt{g}}{\sf SkBrcd}}((G,\cdot,N,\ast,\odot),{\sf R}((H,B,\varphi_{B}))).\end{array}\]

Let $(p,q)\colon {\sf P}((G,\cdot,N,\ast,\odot))=({\sf L}((G,\cdot)),{\sf L}((N,\ast)),\varphi_{\odot})\rightarrow (H,B,\varphi_{B})$ be a morphism in ${}_{\mathbb{K}}\overline{{\sf HBrcd}}$, that is to say, $p\colon \mathbb{K}[G]\rightarrow H$ and $q\colon \mathbb{K}[N]\rightarrow B$ are Hopf algebra morphisms in ${}_{\mathbb{K}}{\sf Vect}$ such that \eqref{Cond. Mor bracoid} holds. Let's define \[{}^{(G,\cdot,N,\ast,\odot)}\Gamma_{(H,B,\varphi_{B})}((p,q))=(\restr{p}{G},\restr{q}{N})\colon (G,\cdot,N,\ast,\odot)\rightarrow {\sf R}((H,B,\varphi_{B}))=({\sf G}(H),{\sf G}(B),\circledast)\]
where $\restr{p}{G}$ and $\restr{q}{N}$ are the restrictions of $p$ and $q$ to $G$ and $N$, respectively. In what follows, we will see that ${}^{(G,\cdot,N,\ast,\odot)}\Gamma_{(H,B,\varphi_{B})}$ is well-defined. In fact, using that both $p$ and $q$ are coalgebra morphisms, we obtain that $p(g)\in{\sf G}(H)$ and $q(n)\in {\sf G}(B)$. Indeed,
\begin{align*}
&\delta_{H}(p(g))\\=&((p\otimes p)\circ\delta_{\mathbb{K}[G]})(g)\;\footnotesize\textnormal{(by the condition of coalgebra morphism for $p$)}\\=&p(g)\otimes p(g),
\end{align*} 
and also,
\begin{align*}
&\varepsilon_{H}(p(g))\\=&\varepsilon_{\mathbb{K}[G]}(g)\;\footnotesize\textnormal{(by the condition of coalgebra morphism for $p$)}\\=&1_{\mathbb{K}},
\end{align*}
so $p(g)$ is a group-like element of $H$. The proof for $q(n)\in{\sf G}(B)$ is analogous. Moreover, due to being $p$ and $q$ algebra morphisms, $\restr{p}{G}\colon G\rightarrow {\sf G}(H)$ and $\restr{q}{N}\colon N\rightarrow {\sf G}(B)$ are group morphisms. We only detail the proof for $\restr{p}{G}$, because for $\restr{q}{N}$ is analogous. Indeed,
 \begin{align*}
 &p(g_{1}\cdot g_{2})=(p\circ\mu_{\mathbb{K}[G]})(g_{1}\otimes g_{2})\\=&(\mu_{H}\circ(p\otimes p))(g_{1}\otimes g_{2})\;\footnotesize\textnormal{(by the condition of algebra morphism for $p$)}\\=&p(g_{1})\diamond p(g_{2}). 
 \end{align*}
 So, to conclude that $(\restr{p}{G},\restr{q}{N})$ is a morphism in $\textnormal{\texttt{g}}\sf{SkBrcd}$, it only remains to check that \eqref{morgenskbracoid} holds. Indeed, consider $g\in G$ and $n\in N$, we have that
 \begin{align*}
 &q(g\odot n)=q(\varphi_{\odot}(g\otimes n))\\=&\varphi_{B}(p(g)\otimes q(n))\;\footnotesize\textnormal{(by \eqref{Cond. Mor bracoid})}\\=&p(g)\circledast q(n).
  \end{align*}
 
 It is a straightforward verification that ${}^{(G,\cdot,N,\ast,\odot)}\Gamma_{(H,B,\varphi_{B})}$ is injective, so to conclude that the correspondence given by ${}^{(G,\cdot,N,\ast,\odot)}\Gamma_{(H,B,\varphi_{B})}$ is bijective let's show that  it is surjective. Consider a morphism $(y,z)\colon (G,\cdot,N,\ast,\odot)\rightarrow {\sf R}((H,B,\varphi_{B}))=({\sf G}(H),{\sf G}(B),\circledast)$ in $\textnormal{\texttt{g}}{\sf SkBrcd}$ and let $\mathbb{K}[y]$ and $\mathbb{K}[z]$ be the linear extensions of $y$ and $z$ to $\mathbb{K}[G]$ and $\mathbb{K}[N]$, respectively. It is easy to check that $\mathbb{K}[y]\colon \mathbb{K}[G]\rightarrow H$ and $\mathbb{K}[z]\colon \mathbb{K}[N]\rightarrow B$ are Hopf algebra morphisms in ${}_{\mathbb{K}}{\sf Vect}$ satisfying \eqref{Cond. Mor bracoid}, condition that follows from \eqref{morgenskbracoid} for $(y,z)$, i.e., $(\mathbb{K}[y],\mathbb{K}[z])$ is a morphism in ${}_{\mathbb{K}}{\sf HBrcd}$. So, the surjectivity of ${}^{(G,\cdot,N,\ast,\odot)}\Gamma_{(H,B,\varphi_{B})}$ follows from:
 \[{}^{(G,\cdot,N,\ast,\odot)}\Gamma_{(H,B,\varphi_{B})}((\mathbb{K}[y],\mathbb{K}[z]))=(\restr{\mathbb{K}[y]}{G},\restr{\mathbb{K}[z]}{N})=(y,z).\]
 
 Therefore, ${\sf P}\dashv{\sf R}$.
\end{proof}
\begin{remark}
Following the arguments and notations of the previous proof, it is straightforward to check that $\left({}^{(G,\cdot,N,\ast,\odot)}\Gamma_{(H,B,\varphi_{B})}\right)^{-1}((y,z))=(\mathbb{K}[y],\mathbb{K}[z])$.
\end{remark}
\begin{definition}
Let $\mathbb{K}$ be a field and $(H,B,\varphi_{B})$ be an object in ${}_{\mathbb{K}}\sf{HBrcd}$. We will say that $(H,B,\varphi_{B})$ is pointed and cosemisimple if the coalgebras $(H,\varepsilon_{H},\delta_{H})$ and $(B,\varepsilon_{B},\delta_{B})$ are pointed and cosemisimple.

Note that by \eqref{pointed+cosemis}, $(H,B,\varphi_{B})$ is pointed and cosemisimple if and only if $H=\mathbb{K}[{\sf G}(H)]$ and $B=\mathbb{K}[{\sf G}(B)]$.
\end{definition}
\begin{corollary}\label{corollaryequiv}
The adjunction of the previous theorem induces a categorical equivalence between $\textnormal{\texttt{g}}{\sf SkBrcd}$ and the full subcategory of ${}_{\mathbb{K}}\overline{{\sf HBrcd}}$ consisting of all pointed and cosemisimple Hopf bracoids.
\end{corollary}
\begin{proof}
We have to see that the unit, $u$, and the counit, $v$, of the adjunction defined in Theorem \ref{adjunction} are natural isomorphisms. 

On the one hand, for an object $(G,\cdot,N,\ast,\odot)$ in  $\textnormal{\texttt{g}}{\sf SkBrcd}$, the unit is given by
\[u_{(G,\cdot,N,\ast,\odot)}={}^{(G,\cdot,N,\ast,\odot)}\Gamma_{{\sf P}((G,\cdot,N,\ast,\odot))}(id_{{\sf P}((G,\cdot,N,\ast,\odot))})\colon (G,\cdot,N,\ast,\odot)\rightarrow ({\sf R}\circ{\sf P})(G,\cdot,N,\ast,\odot).\] 
Note that 
\[({\sf R}\circ{\sf P})(G,\cdot,N,\ast,\odot)={\sf R}((\mathbb{K}[G],\mathbb{K}[N],\varphi_{\odot}))=({\sf G}(\mathbb{K}[G]),{\sf G}(\mathbb{K}[N]),\circledast)=(G,N,\circledast)\]
where, in this case, $g\circledast n=\varphi_{\odot}(g\otimes n)=g\odot n$ for all $g\in G$ and $n \in N$, that is to say, $\circledast=\odot$, and therefore, ${\sf R}\circ{\sf P}={\sf id}_{\textnormal{\texttt{g}}{\sf SkBrcd}}$. Moreover, 
\begin{align*}
&{}^{(G,\cdot,N,\ast,\odot)}\Gamma_{{\sf P}((G,\cdot,N,\ast,\odot))}(id_{{\sf P}((G,\cdot,N,\ast,\odot))})={}^{(G,\cdot,N,\ast,\odot)}\Gamma_{(\mathbb{K}[G],\mathbb{K}[N],\varphi_{\odot})}((id_{\mathbb{K}[G]},id_{\mathbb{K}[N]}))\\=&(\restr{id_{\mathbb{K}[G]}}{G},\restr{id_{\mathbb{K}[N]}}{N})=(id_{G},id_{N}).
\end{align*}
So, as a conclusion, $u_{(G,\cdot,N,\ast,\odot)}=(id_{G},id_{N})=id_{(G,\cdot,N,\ast,\odot)}$, which is obvious an isomorphism in $\texttt{g}{\sf SkBrcd}$.

On the other hand, the counit is defined as follows for a pointed and cosemisimple Hopf bracoid $(H,B,\varphi_{B})$ in ${}_{\mathbb{K}}\overline{{\sf HBrcd}}$:
\[v_{(H,B,\varphi_{B})}=\left({}^{{\sf R}((H,B,\varphi_{B}))}\Gamma_{(H,B,\varphi_{B})}\right)^{-1}(id_{{\sf R}((H,B,\varphi_{B}))})\colon ({\sf P}\circ{\sf R})(H,B,\varphi_{B})\rightarrow (H,B,\varphi_{B}).\]
Thanks to the pointed and cosemisimple character, we obtain that
\[({\sf P}\circ{\sf R})(H,B,\varphi_{B})={\sf P}(({\sf G}(H),{\sf G}(B),\circledast))=(\mathbb{K}[{\sf G}(H)],\mathbb{K}[{\sf G}(B)],\varphi_{\circledast})=(H,B,\varphi_{\circledast})\]
where $\varphi_{\circledast}(h\otimes b)=h\circledast b=\varphi_{B}(h\otimes b)$ for all $h\in {\sf G}(H)$ and $b\in{\sf G}(B)$, which implies that $\varphi_{\circledast}=\varphi_{B}$, and thus $({\sf P}\circ{\sf R})(H,B,\varphi_{B})=(H,B,\varphi_{B})$. In addition, 
\begin{align*}
&\left({}^{{\sf R}((H,B,\varphi_{B}))}\Gamma_{(H,B,\varphi_{B})}\right)^{-1}(id_{{\sf R}((H,B,\varphi_{B}))})=\left({}^{({\sf G}(H),{\sf G}(B),\circledast)}\Gamma_{(H,B,\varphi_{B})}\right)^{-1}((id_{{\sf G}(H)}, id_{{\sf G}(B)}))\\=&(\mathbb{K}[id_{{\sf G}(H)}],\mathbb{K}[id_{{\sf G}(B)}])=(id_{\mathbb{K}[{\sf G}(H)]},id_{\mathbb{K}[{\sf G}(B)]})\\=&(id_{H},id_{B})\;\footnotesize\textnormal{(by the pointed and consemisimple character of $(H,B,\varphi_{B})$)}.
\end{align*}
Therefore, $v_{(H,B,\varphi_{B})}=(id_{H},id_{B})=id_{(H,B,\varphi_{B})}$, which is obvious an isomorphism in ${}_{\mathbb{K}}{\sf HBrcd}$.

Due to being the unit and the counit of the adjunction natural isomorphisms, we obtain the required categorical equivalence.
\end{proof}

\begin{lemma}\label{opbracoid}
Let's assume that $\sf{C}$ is symmetric. Let $(H,B,\varphi_{B})$ be a Hopf bracoid in $\sf{C}$ with $H$ cocommutative. If $\overline{B}=(B,\eta_{B},\mu_{\overline{B}},\varepsilon_{B},\delta_{B},\lambda_{B})$ is a Hopf algebra in $\sf{C}$, then $(H,\overline{B},\varphi_{B})$ is also a Hopf bracoid in $\sf{C}$, being $\mu_{\overline{B}}\coloneqq \mu_{B}\circ c_{B,B}$.
\end{lemma}
\begin{proof}
It is enough to see that \eqref{Eq. Bracoid compatibility} holds for $(H,\overline{B},\varphi_{B})$. Indeed:
\begin{align*}
&\mu_{\overline{B}}\circ(\varphi_{B}\otimes\Gamma_{\overline{B}})\circ(H\otimes c_{H,B}\otimes B)\circ(\delta_{H}\otimes B\otimes B)\\=&\mu_{B}\circ(\varphi_{B}\otimes (\mu_{B}\circ((\lambda_{B}\circ u_{\varphi_{B}})\otimes\varphi_{B})))\circ(H\otimes((c_{H,B}\otimes H)\circ(H\otimes c_{H,B})\circ(\delta_{H}\otimes B))\otimes B)\\&\circ(\delta_{H}\otimes c_{B,B})\;\footnotesize\textnormal{(by naturality of $c$, symmetry and cocommutativity of $\delta_{H}$)}\\=&\mu_{B}\circ(\varphi_{B}\otimes \Phi_{B})\circ(H\otimes c_{H,B}\otimes B)\circ(\delta_{H}\otimes c_{B,B})\;\footnotesize\textnormal{(by naturality of $c$ and definition of $\Phi_{B}$)}\\=&\varphi_{B}\circ(H\otimes(\mu_{B}\circ c_{B,B}))\;\footnotesize\textnormal{(by \eqref{Eq. Bracoid compatibility} for $(H,B,\varphi_{B})$)}\\=&\varphi_{B}\circ(H\otimes\mu_{\overline{B}}).
\qedhere
\end{align*}

\end{proof}
\begin{remark}\label{overB-HA}
Note that $\overline{B}=(B,\eta_{B},\mu_{\overline{B}},\delta_{B},\varepsilon_{B},\lambda_{B})$ is not always a Hopf algebra in $\sf{C}$, but it does under cocommutativity of $B$ or when $\lambda_{B}\circ\lambda_{B}=id_{B}$. 
\end{remark}
\begin{example}
Take ${\sf C}={\sf Set}$, which is a symmetric monoidal category whose symmetry $c$ is the usual flip, $c(x,y)=(y,x)$ for all $x\in X$ and $y\in Y$. By Example \ref{example functor P}, a Hopf bracoid $(H,B,\varphi_{B})$ in ${\sf Set}$ is no more than a generalized skew bracoid, and also every Hopf algebra $H$ in ${\sf Set}$ is cocommutative because the coproduct of $H$ is given by duplication, $\delta_{H}(h)=(h,h)$. Moreover, the antipode of the Hopf algebra $B$ in ${\sf Set}$ is the inversion, $\lambda_{B}(b)=b^{-1}$, so $\lambda_{B}\circ\lambda_{B}=id_{B}$, and then $\overline{B}$ is also a Hopf algebra in ${\sf Set}$ by Remark \ref{overB-HA}. Therefore, by Lemma \ref{opbracoid}, if $(H,B,\varphi_{B})$ is a Hopf bracoid in {\sf Set}, then $(H,\overline{B},\varphi_{B})$ is a Hopf bracoid in {\sf Set}, which is exactly the result obtained by Martin-Lyons and Truman in \cite[Proposition 2.7]{MLT}.
\end{example}
\begin{lemma}\label{Despejar phiB}
If $(H,B,\varphi_{B})$ is a Hopf bracoid in $\sf{C}$, then the equality
\begin{equation}\label{phiB despejar}
\varphi_{B}=\mu_{B}\circ(u_{\varphi_{B}}\otimes\Phi_{B})\circ(\delta_{H}\otimes B).
\end{equation}
holds.
\end{lemma}
\begin{proof}
Composing in \eqref{Eq. Bracoid compatibility} with $H\otimes\eta_{B}\otimes B$ on the right and using the unit property $\mu_{B}\circ(\eta_{B}\otimes B)=id_{B}$ and the naturality of $c$ we obtain the required result.
\end{proof}
\begin{remark}\label{rmkcoalgu1}
Note that 
\begin{align}\label{rmk1}
&\mu_{B}\circ(u_{\varphi_{B}}\otimes\Phi_{B})\circ(\delta_{H}\otimes B)\\\nonumber=&\mu_{B}\circ((u_{\varphi_{B}}\ast(\lambda_{B}\circ u_{\varphi_{B}}))\otimes\varphi_{B})\circ(\delta_{H}\otimes B)\;\footnotesize\textnormal{(by coassociativity of $\delta_{H}$ and associativity of $\mu_{B}$)}.
\end{align}
If we suppose that $\varphi_{B}$ is a coalgebra morphism, then $u_{\varphi_{B}}$ is a coalgebra morphism, that is to say, the equalities 
\begin{align}
\label{u1coalgdelta} (u_{\varphi_{B}}\otimes u_{\varphi_{B}})\circ\delta_{H}=\delta_{B}\circ u_{\varphi_{B}},\\
\label{u1coalgeps} \varepsilon_{B}\circ u_{\varphi_{B}}=\varepsilon_{H}
\end{align}
hold. Indeed, the proof of \eqref{u1coalgeps} is straightforward and \eqref{u1coalgdelta} follows by:
\begin{align*}
&(u_{\varphi_{B}}\otimes u_{\varphi_{B}})\circ\delta_{H}\\=&(\varphi_{B}\otimes\varphi_{B})\circ(H\otimes c_{H,B}\otimes B)\circ(\delta_{H}\otimes\eta_{B}\otimes\eta_{B})\;\footnotesize\textnormal{(by definition of $u_{\varphi_{B}}$ and naturality of $c$)}\\=&(\varphi_{B}\otimes\varphi_{B})\circ(H\otimes c_{H,B}\otimes B)\circ(\delta_{H}\otimes(\delta_{B}\circ\eta_{B}))\;\footnotesize\textnormal{(by the condition of coalgebra morphism for $\eta_{B}$)}\\=&\delta_{B}\circ\varphi_{B}\circ(H\otimes\eta_{B})\;\footnotesize\textnormal{(by the condition of coalgebra morphism for $\varphi_{B}$)}\\=&\delta_{B}\circ u_{\varphi_{B}}.
\end{align*}
Therefore, by the condition of coalgebra morphism for $u_{\varphi_{B}}$, we obtain that
\begin{align}\label{rmk2}
&u_{\varphi_{B}}\ast(\lambda_{B}\circ u_{\varphi_{B}})\\\nonumber=&(id_{B}\ast\lambda_{B})\circ u_{\varphi_{B}}\;\footnotesize\textnormal{(by \eqref{u1coalgdelta})}\\\nonumber=&\eta_{B}\circ\varepsilon_{B}\circ u_{\varphi_{B}}\;\footnotesize\textnormal{(by \eqref{antipode})}\\\nonumber=&\varepsilon_{H}\otimes\eta_{B}\;\footnotesize\textnormal{(by \eqref{u1coalgeps})}.
\end{align}
So, we conclude that
\begin{align*}
&\mu_{B}\circ(u_{\varphi_{B}}\otimes\Phi_{B})\circ(\delta_{H}\otimes B)\\=&\mu_{B}\circ((u_{\varphi_{B}}\ast(\lambda_{B}\circ u_{\varphi_{B}}))\otimes\varphi_{B})\circ(\delta_{H}\otimes B)\;\footnotesize\textnormal{(by \eqref{rmk1})}\\=&\mu_{B}\circ((\eta_{B}\circ\varepsilon_{H})\otimes\varphi_{B})\circ(\delta_{H}\otimes B)\;\footnotesize\textnormal{(by \eqref{rmk2})}\\=&\varphi_{B}\;\footnotesize\textnormal{(by (co)unit properties)},
\end{align*}
i.e., when $\varphi_{B}$ is a coalgebra morphism, \eqref{phiB despejar} always holds, without the need for the condition \eqref{Eq. Bracoid compatibility} to be fulfilled. 
\end{remark}

From now on, we are going to denote by $\overline{\sf{HBrcd}}$ to the full subcategory of $\sf{HBrcd}$ whose objects are Hopf bracoids $(H,B,\varphi_{B})$ such that $\varphi_{B}$ is a coalgebra morphism, i.e., $(B,\varphi_{B})$ is a left $H$-module coalgebra.
%\color{red}
%($\varphi_{B}$ no es posible despejarlo en función de $\Phi_{B}$ y $\mu_{B}$, sino que aparece $\varphi_{B}$.)
%\color{black}

It is well known that $(H_{1},\Gamma_{H_{1}})$ is a left $H_{2}$-module algebra in $\sf{C}$ when $\mathbb{H}=(H_{1},H_{2})$ is a Hopf brace in $\sf{C}$. So, it is natural to wonder what conditions does $\Phi_{B}$ satisfy when $(H,B,\varphi_{B})$ is a Hopf bracoid in $\sf{C}$. First of all, the following lemma will be relevant in the subsequent results, which is the generalization of \cite[Lemma 1.7]{AGV} in the context of Hopf bracoids.
\begin{lemma}
Let $(H,B,\varphi_{B})$ be a Hopf bracoid in $\sf{C}$. If $\varphi_{B}$ is a coalgebra morphism, then 
\begin{equation}\label{Cond Lem Imp}
\Phi_{B}\circ(H\otimes\lambda_{B})=\mu_{B}\circ((\lambda_{B}\circ \varphi_{B})\otimes u_{\varphi_{B}})\circ(H\otimes c_{H,B})\circ(\delta_{H}\otimes B).
\end{equation}
\end{lemma}
\begin{proof} By \eqref{antipode} for $B$ and \eqref{Eq. Bracoid compatibility} we have that
\begin{equation}\label{lem1}
\varphi_{B}\circ(H\otimes(\eta_{B}\circ\varepsilon_{B}))=\mu_{B}\circ(\varphi_{B}\otimes\Phi_{B})\circ(H\otimes c_{H,B}\otimes \lambda_{B})\circ(\delta_{H}\otimes\delta_{B}).
\end{equation}
So, using \eqref{lem1}, we obtain the following:
\begin{align*}
&\Phi_{B}\circ(H\otimes\lambda_{B})\\=&((\varepsilon_{B}\circ\varphi_{B})\otimes\Phi_{B})\circ(H\otimes c_{H,B}\otimes\lambda_{B})\circ(\delta_{H}\otimes\delta_{B})\;\footnotesize\textnormal{(by the condition of coalgebra morphism for $\varphi_{B}$,}\\&\footnotesize\textnormal{naturality of $c$ and counit properties)}\\=&\mu_{B}\circ((\eta_{B}\circ\varepsilon_{B}\circ\varphi_{B})\otimes\Phi_{B})\circ(H\otimes c_{H,B}\otimes\lambda_{B})\circ(\delta_{H}\otimes\delta_{B})\;\footnotesize\textnormal{(by unit property)}\\=&\mu_{B}\circ((\mu_{B}\circ(\lambda_{B}\otimes B)\circ\delta_{B}\circ\varphi_{B})\otimes\Phi_{B})\circ(H\otimes c_{H,B}\otimes\lambda_{B})\circ(\delta_{H}\otimes\delta_{B})\;\footnotesize\textnormal{(by \eqref{antipode})}\\=&\mu_{B}\circ((\mu_{B}\circ(\lambda_{B}\otimes B)\circ(\varphi_{B}\otimes\varphi_{B})\circ(H\otimes c_{H,B}\otimes B)\circ(\delta_{H}\otimes\delta_{B}))\otimes\Phi_{B})\circ(H\otimes c_{H,B}\otimes \lambda_{B})\\&\circ(\delta_{H}\otimes\delta_{B})\;\footnotesize\textnormal{(by the condition of coalgebra morphism for $\varphi_{B}$)}\\=&\mu_{B}\circ((\lambda_{B}\circ \varphi_{B})\otimes(\mu_{B}\circ(\varphi_{B}\otimes\Phi_{B})\circ(H\otimes c_{H,B}\otimes B)))\circ (H\otimes((c_{H,B}\otimes H)\circ(H\otimes c_{H,B})\circ(\delta_{H}\otimes B))\\&\otimes((B\otimes\lambda_{B})\circ\delta_{B}))\circ(\delta_{H}\otimes\delta_{B})\;\footnotesize\textnormal{(by naturality of $c$, associativity of $\mu_{B}$ and coassociativity of $\delta_{H}$ and $\delta_{B}$)}\\=&\mu_{B}\circ((\lambda_{B}\circ\varphi_{B})\otimes(\mu_{B}\circ(\varphi_{B}\otimes\Phi_{B})\circ(H\otimes c_{H,B}\otimes \lambda_{B})\circ(\delta_{H}\otimes \delta_{B})))\circ(H\otimes c_{H,B}\otimes B)\circ(\delta_{H}\otimes\delta_{B})\\&\footnotesize\textnormal{(by naturality of $c$)}\\=&\mu_{B}\circ((\lambda_{B}\circ\varphi_{B})\otimes(\varphi_{B}\circ(H\otimes(\eta_{B}\circ\varepsilon_{B}))))\circ(H\otimes c_{H,B}\otimes B)\circ(\delta_{H}\otimes\delta_{B})\;\footnotesize\textnormal{(by \eqref{lem1})}\\=&\mu_{B}\circ((\lambda_{B}\circ\varphi_{B})\otimes u_{\varphi_{B}})\circ(H\otimes c_{H,B})\circ(\delta_{H}\otimes B)\;\footnotesize\textnormal{(by counit property and definition of $u_{\varphi_{B}}$)}.\qedhere
\end{align*}
\end{proof}
%\color{red}
%A lo mejor, como la condición ``$\varphi_{B}$ morfismo de coálgebras'' va a ser recurrente en lo sucesivo, deberíamos plantearnos incluirla en la definición de Hopf bracoid.
%\color{black}
\begin{theorem}\label{GammaB modulo algebra}
Let $(H,B,\varphi_{B})$ be a Hopf bracoid such that $\varphi_{B}$ is a coalgebra morphism, then $(B,\Phi_{B})$ is a left $H$-module algebra.
\end{theorem}
\begin{proof}
First of all, let's see that $\Phi_{B}\circ(\eta_{H}\otimes B)=id_{B}$. Note that, by the $H$-module conditions, it is straightforward to prove that
\begin{align}\label{etau1}
u_{\varphi_{B}}\circ\eta_{H}=\eta_{B},\\\label{phiBu1}\varphi_{B}\circ(H\otimes u_{\varphi_{B}})=u_{\varphi_{B}}\circ\mu_{H}.
\end{align}
Then,
\begin{align*}
&\Phi_{B}\circ(\eta_{H}\otimes B)\\=&\mu_{B}\circ((\lambda_{B}\circ u_{\varphi_{B}}\circ\eta_{H})\otimes(\varphi_{B}\circ(\eta_{H}\otimes B)))\;\footnotesize\textnormal{(by the condition of coalgebra morphism for $\eta_{H}$)}\\=&id_{B}\;\footnotesize\textnormal{(by \eqref{etau1}, \eqref{u-antip}, $H$-module conditions and unit property)}.
\end{align*}
Moreover, $\Phi_{B}\circ(H\otimes \Phi_{B})=\Phi_{B}\circ(\mu_{H}\otimes B)$ follows by:
\begin{align*}
&\Phi_{B}\circ(H\otimes \Phi_{B})\\=&\mu_{B}\circ((\Phi_{B}\circ(H\otimes\lambda_{B}))\otimes\Phi_{B})\circ(H\otimes c_{H,B}\otimes B)\circ(\delta_{H}\otimes((u_{\varphi_{B}}\otimes\varphi_{B})\circ(\delta_{H}\otimes B)))\;\footnotesize\textnormal{(by \eqref{Eq. Bracoid compatibility},}\\&\footnotesize\textnormal{coassociativity of $\delta_{H}$, definition of $\Phi_{B}$ and naturality of $c$)}\\=&\mu_{B}\circ((\mu_{B}\circ((\lambda_{B}\circ\varphi_{B})\otimes u_{\varphi_{B}}))\otimes\Phi_{B})\circ(H\otimes((c_{H,B}\otimes H)\circ(H\otimes c_{H,B})\circ(\delta_{H}\otimes B))\otimes B)\\&\circ(\delta_{H}\otimes((u_{\varphi_{B}}\otimes\varphi_{B})\circ(\delta_{H}\otimes B)))\;\footnotesize\textnormal{(by \eqref{Cond Lem Imp} and coassociativity of $\delta_{H}$)}\\=&\mu_{B}\circ((\lambda_{B}\circ\varphi_{B})\otimes(\mu_{B}\circ (u_{\varphi_{B}}\otimes\Phi_{B})\circ(\delta_{H}\otimes B)))\circ(H\otimes c_{H,B}\otimes B)\circ(\delta_{H}\otimes((u_{\varphi_{B}}\otimes\varphi_{B})\\& \circ(\delta_{H}\otimes B)))\footnotesize\textnormal{(by naturality of $c$ and associativity of $\mu_{B}$)}\\=&\mu_{B}\circ((\lambda_{B}\circ\varphi_{B}\circ(H\otimes u_{\varphi_{B}}))\otimes(\varphi_{B}\circ(H\otimes\varphi_{B})))\circ(H\otimes c_{H,H}\otimes H\otimes B)\circ(\delta_{H}\otimes\delta_{H}\otimes B)\;\footnotesize\textnormal{(by \eqref{phiB despejar}}\\&\footnotesize\textnormal{and naturality of $c$)}\\=&\mu_{B}\circ((\lambda_{B}\circ u_{\varphi_{B}})\otimes\varphi_{B})\circ(((\mu_{H}\otimes\mu_{H})\circ(H\otimes c_{H,H}\otimes H)\circ(\delta_{H}\otimes\delta_{H}))\otimes B)\;\footnotesize\textnormal{(by $H$-module conditions}\\&\footnotesize\textnormal{and \eqref{phiBu1})}\\=&\Phi_{B}\circ(\mu_{H}\otimes B)\;\footnotesize\textnormal{(by the condition of coalgebra morphsim for $\mu_{H}$ and definition of $\Phi_{B}$)}.
\end{align*}
So we conclude that $(B,\Phi_{B})$ is a left $H$-module. To finish the proof, we have to compute that $\mu_{B}$ and $\eta_{B}$ are morphism of left $H$-modules. On the one side,
\begin{align*}
&\Phi_{B}\circ(H\otimes\eta_{B})=\mu_{B}\circ((\lambda_{B}\circ u_{\varphi_{B}})\otimes u_{\varphi_{B}})\circ\delta_{H}\\=&(\lambda_{B}\ast id_{B})\circ u_{\varphi_{B}}\;\footnotesize\textnormal{(by \eqref{u1coalgdelta})}\\=&\varepsilon_{H}\otimes\eta_{B}\;\footnotesize\textnormal{(by \eqref{antipode} and \eqref{u1coalgeps})}.
\end{align*}
On the other side,
\begin{align*}
&\Phi_{B}\circ(H\otimes\mu_{B})\\=&\mu_{B}\circ((\lambda_{B}\circ u_{\varphi_{B}})\otimes(\mu_{B}\circ(\varphi_{B}\otimes\Phi_{B})\circ(H\otimes c_{H,B}\otimes B)\circ(\delta_{H}\otimes B\otimes B)))\circ(\delta_{H}\otimes B\otimes B)\;\footnotesize\textnormal{(by \eqref{Eq. Bracoid compatibility})}\\=&\mu_{B}\circ((\mu_{B}\circ((\lambda_{B}\circ u_{\varphi_{B}})\otimes\varphi_{B})\circ(\delta_{H}\otimes B))\otimes\Phi_{B})\circ(H\otimes c_{H,B}\otimes B)\circ(\delta_{H}\otimes B\otimes B)\;\footnotesize\textnormal{(by coassociativity}\\&\footnotesize\textnormal{of $\delta_{H}$ and associativity of $\mu_{B}$)}\\=&\mu_{B}\circ(\Phi_{B}\otimes\Phi_{B})\circ(H\otimes c_{H,B}\otimes B)\circ(\delta_{H}\otimes B\otimes B).
\qedhere
\end{align*}
\end{proof}
\begin{lemma}\label{GammaBcoalg}
Let $(H,B,\varphi_{B})$ be a Hopf bracoid in $\sf{C}$ such that $\varphi_{B}$ is a coalgebra morphism. If 
\begin{equation}\label{clas coc}
(\Phi_{B}\otimes H)\circ(H\otimes c_{H,B})\circ((c_{H,H}\circ \delta_{H})\otimes B)=(\Phi_{B}\otimes H)\circ(H\otimes c_{H,B})\circ(\delta_{H}\otimes B),
\end{equation}
then $\Phi_{B}$ is a coalgebra morphism.
\end{lemma}
\begin{proof}
Note that $\varepsilon_{B}\circ\Phi_{B}=\varepsilon_{H}\otimes\varepsilon_{B}$ holds thanks to the condition of coalgebra morphism for $\mu_{B}$, \eqref{u-antip}, \eqref{u1coalgeps} and the condition of coalgebra morphism for $\varphi_{B}$. Moreover,
\begin{align*}
&\delta_{B}\circ\Phi_{B}\\=&(\mu_{B}\otimes\mu_{B})\circ(B\otimes c_{B,B}\otimes B)\circ((\delta_{B}\circ\lambda_{B}\circ u_{\varphi_{B}})\otimes(\delta_{B}\circ\varphi_{B}))\circ(\delta_{H}\otimes B)\;\footnotesize\textnormal{(by the condition of coalgebra}\\&\footnotesize\textnormal{morphism for $\mu_{B}$)}\\=&(\mu_{B}\otimes\mu_{B})\circ(B\otimes c_{B,B}\otimes B)\circ((c_{B,B}\circ(\lambda_{B}\otimes\lambda_{B})\circ(u_{\varphi_{B}}\otimes u_{\varphi_{B}})\circ\delta_{H})\otimes((\varphi_{B}\otimes\varphi_{B})\circ(H\otimes c_{H,B}\otimes B)\\&\circ(\delta_{H}\otimes\delta_{B})))\circ(\delta_{H}\otimes B)\;\footnotesize\textnormal{(by the condition of coalgebra morphism for $\varphi_{B}$, \eqref{a-antip} and \eqref{u1coalgdelta})}\\=&((\mu_{B}\circ((\lambda_{B}\circ u_{\varphi_{B}})\otimes \varphi_{B}))\otimes(\mu_{B}\circ((\lambda_{B}\circ u_{\varphi_{B}})\otimes \varphi_{B})))\circ(H\otimes H\otimes c_{H,B}\otimes H\otimes B)\\&\circ(((H\otimes c_{H,H})\circ(c_{H,H}\otimes H)\circ(H\otimes\delta_{H}))\otimes c_{H,B}\otimes B)\circ(((\delta_{H}\otimes H)\circ\delta_{H})\otimes\delta_{B})\;\footnotesize\textnormal{(by naturality of $c$ and}\\&\footnotesize\textnormal{coassociativity of $\delta_{H}$)}\\=&(B\otimes(\mu_{B}\circ((\lambda_{B}\circ u_{\varphi_{B}})\otimes\varphi_{B})))\circ(((\Phi_{B}\otimes H)\circ(H\otimes c_{H,B})\circ((c_{H,H}\circ\delta_{H})\otimes B))\otimes H\otimes B)\\&\circ(H\otimes c_{H,B}\otimes B)\circ(\delta_{H}\otimes\delta_{B})\;\footnotesize\textnormal{(by naturality of $c$ and definition of $\Phi_{B}$)}\\=&(\Phi_{B}\otimes(\mu_{B}\circ((\lambda_{B}\circ u_{\varphi_{B}})\otimes\varphi_{B})))\circ(H\otimes((c_{H,B}\otimes H)\circ(H\otimes c_{H,B})\circ(\delta_{H}\otimes B))\otimes B)\circ(\delta_{H}\otimes\delta_{B})\\&\footnotesize\textnormal{(by \eqref{clas coc} and coassociativity of $\delta_{H}$)}\\=&(\Phi_{B}\otimes\Phi_{B})\circ(H\otimes c_{H,B}\otimes B)\circ(\delta_{H}\otimes\delta_{B})\;\footnotesize\textnormal{(by naturality of $c$ and definition of $\Phi_{B}$)}.
\qedhere
\end{align*}
\begin{remark}\label{cc1}
If {\sf C} is symmetric, \eqref{clas coc} means that $(B,\Phi_{B})$ is in the cocommutativity class of $H$ following \cite[Definition 2.1 and 2.2]{CCH}. On the other hand, note that when $H$ is a cocommutative Hopf algebra, then \eqref{clas coc} always holds.
\end{remark}
\end{proof}
Note that we can rewrite right side of condition \eqref{Eq. Bracoid compatibility} as follows:
%\footnotesize
\begin{align*}
&\mu_{B}\circ(\varphi_{B}\otimes\Phi_{B})\circ(H\otimes c_{H,B}\otimes B)\circ(\delta_{H}\otimes B\otimes B)\\=&\mu_{B}\circ((\mu_{B}\circ(\varphi_{B}\otimes(\lambda_{B}\circ u_{\varphi_{B}}))\circ (H\otimes c_{H,B})\circ(\delta_{H}\otimes B))\otimes\varphi_{B})\circ(H\otimes c_{H,B}\otimes B)\circ(\delta_{H}\otimes B\otimes B)\\&\footnotesize\textnormal{(by naturality of $c$, coassociativity of $\delta_{H}$ and associativity of $\mu_{B}$)}\\=&\mu_{B}\circ(\Phi'_{B}\otimes\varphi_{B})\circ(H\otimes c_{H,B}\otimes B)\circ(\delta_{H}\otimes B\otimes B)
\end{align*}
%\normalsize
where $\Phi'_{B}\coloneqq \mu_{B}\circ(\varphi_{B}\otimes(\lambda_{B}\circ u_{\varphi_{B}}))\circ(H\otimes c_{H,B})\circ(\delta_{H}\otimes B)$. Therefore, compatibility condition \eqref{Eq. Bracoid compatibility} of Definition \ref{HBrcd} is equivalente to:
\begin{equation}\label{Compat prima}
\varphi_{B}\circ(H\otimes\mu_{B})=\mu_{B}\circ(\Phi'_{B}\otimes\varphi_{B})\circ(H\otimes c_{H,B}\otimes B)\circ(\delta_{H}\otimes B\otimes B).
\end{equation}

In the following results we will see that $(B,\Phi'_{B})$ is also a left $H$-module algebra under certain conditions. Firstly, take into account the following lemma:
\begin{lemma}
Let $(H,B,\varphi_{B})$ be a Hopf bracoid in $\sf{C}$ such that $\varphi_{B}$ is a coalgebra morphism, then
\begin{equation}\label{imp for gammaBprima}
\lambda_{B}\circ\varphi_{B}\circ(H\otimes u_{\varphi_{B}})=\mu_{B}\circ(\Phi_{B}\otimes B)\circ(H\otimes((\lambda_{B}\otimes\lambda_{B})\circ(u_{\varphi_{B}}\otimes u_{\varphi_{B}})\circ c_{H,H}))\circ(\delta_{H}\otimes H).
\end{equation}
\begin{proof}
The proof follows by:
\begin{align*}
&\mu_{B}\circ(\Phi_{B}\otimes B)\circ(H\otimes((\lambda_{B}\otimes\lambda_{B})\circ(u_{\varphi_{B}}\otimes u_{\varphi_{B}})\circ c_{H,H}))\circ(\delta_{H}\otimes H)\\=&\mu_{B}\circ((\mu_{B}\circ ((\lambda_{B}\circ\varphi_{B}\circ(H\otimes u_{\varphi_{B}}))\otimes u_{\varphi_{B}})\circ(H\otimes c_{H,H})\circ(\delta_{H}\otimes H))\otimes (\lambda_{B}\circ u_{\varphi_{B}}))\\&\circ(H\otimes c_{H,H})\circ(\delta_{H}\otimes H)\;\footnotesize\textnormal{(by \eqref{Cond Lem Imp} and naturality of $c$)}\\=&\mu_{B}\circ((\lambda_{B}\circ\varphi_{B}\circ(H\otimes u_{\varphi_{B}}))\otimes(\mu_{B}\circ (u_{\varphi_{B}}\otimes(\lambda_{B}\circ u_{\varphi_{B}}))))\circ(H\otimes ((c_{H,H}\otimes H)\circ (H\otimes c_{H,H})\\&\circ(\delta_{H}\otimes H)))\circ(\delta_{H}\otimes H)\;\footnotesize\textnormal{(by associativity of $\mu_{B}$ and coassociativity of $\delta_{H}$)}\\=&\mu_{B}\circ((\lambda_{B}\circ\varphi_{B}\circ(H\otimes u_{\varphi_{B}}))\otimes(u_{\varphi_{B}}\ast(\lambda_{B}\circ u_{\varphi_{B}})))\circ(H\otimes c_{H,H})\circ(\delta_{H}\otimes H)\;\footnotesize\textnormal{(by naturality of $c$)}\\=&\mu_{B}\circ((\lambda_{B}\circ\varphi_{B}\circ(H\otimes u_{\varphi_{B}}))\otimes(\eta_{B}\circ\varepsilon_{H}))\circ(H\otimes c_{H,H})\circ(\delta_{H}\otimes H)\;\footnotesize\textnormal{(by \eqref{rmk2})}\\=&\lambda_{B}\circ\varphi_{B}\circ(H\otimes u_{\varphi_{B}})\;\footnotesize\textnormal{(by (co)unit properties and naturality of $c$)}.\qedhere
\end{align*}
\end{proof}
\end{lemma}
\begin{theorem}
Let $(H,B,\varphi_{B})$ be a Hopf bracoid in $\sf{C}$ such that $\varphi_{B}$ is a coalgebra morphism, then $(B,\Phi'_{B})$ is a left $H$-module algebra.
\end{theorem}
\begin{proof}
At first, we will see that $\Phi'_{B}\circ(\eta_{H}\otimes B)=id_{B}$.
\begin{align*}
&\Phi'_{B}\circ(\eta_{H}\otimes B)\\=&\mu_{B}\circ((\varphi_{B}\circ (\eta_{H}\otimes B))\otimes(\lambda_{B}\circ u_{\varphi_{B}}\circ\eta_{H}))\;\footnotesize\textnormal{(by the condition of coalgebra morphism for $\eta_{H}$}\\&\footnotesize\textnormal{and naturality of $c$)}\\=&id_{B}\;\footnotesize\textnormal{(by the conditions of left $H$-module, \eqref{etau1}, \eqref{u-antip} and unit property)}.
\end{align*}
In addition, the equality $\Phi'_{B}\circ(H\otimes\Phi'_{B})=\Phi'_{B}\circ(\mu_{H}\otimes B)$ also holds. Indeed:
\begin{align*}
&\Phi'_{B}\circ(H\otimes\Phi'_{B})\\=&\mu_{B}\circ((\varphi_{B}\circ(H\otimes\mu_{B}))\otimes(\lambda_{B}\circ u_{\varphi_{B}}))\circ(H\otimes((\varphi_{B}\otimes c_{H,B})\circ(H\otimes c_{H,B}\otimes(\lambda_{B}\circ u_{\varphi_{B}}))\\&\circ(c_{H,H}\otimes c_{H,B})))\circ(\delta_{H}\otimes\delta_{H}\otimes B)\;\footnotesize\textnormal{(by definition of $\Phi'_{B}$ and naturality of $c$)}\\=&\mu_{B}\circ((\mu_{B}\circ(\varphi_{B}\otimes\Phi_{B})\circ(H\otimes c_{H,B}\otimes B)\circ(\delta_{H}\otimes B\otimes B))\otimes(\lambda_{B}\circ u_{\varphi_{B}}))\circ(H\otimes((\varphi_{B}\otimes c_{H,B})\\&\circ(H\otimes c_{H,B}\otimes(\lambda_{B}\circ u_{\varphi_{B}}))\circ(c_{H,H}\otimes c_{H,B})))\circ(\delta_{H}\otimes\delta_{H}\otimes B)\;\footnotesize\textnormal{(by \eqref{Eq. Bracoid compatibility})}\\=&\mu_{B}\circ((\varphi_{B}\circ(\mu_{H}\otimes B))\otimes(\mu_{B}\circ(\Phi_{B}\otimes B)\circ(H\otimes((\lambda_{B}\otimes \lambda_{B})\circ (u_{\varphi_{B}}\otimes u_{\varphi_{B}})\circ c_{H,H}))\circ(\delta_{H}\otimes H)))\\&\circ(H\otimes((H\otimes c_{H,B}\otimes H)\circ(c_{H,H}\otimes c_{H,B})))\circ(\delta_{H}\otimes\delta_{H}\otimes B)\;\footnotesize\textnormal{(by associativity of $\mu_{B}$, coassociativity of $\delta_{H}$,}\\&\footnotesize\textnormal{naturality of $c$ and the conditions of left $H$-module)}\\=&\mu_{B}\circ((\varphi_{B}\circ(\mu_{H}\otimes B))\otimes(\lambda_{B}\circ u_{\varphi_{B}}\circ\mu_{H}))\circ(H\otimes((H\otimes c_{H,B}\otimes H)\circ(c_{H,H}\otimes c_{H,B})))\\&\circ(\delta_{H}\otimes\delta_{H}\otimes B)\;\footnotesize\textnormal{(by \eqref{imp for gammaBprima} and \eqref{phiBu1})}\\=&\mu_{B}\circ(\varphi_{B}\otimes(\lambda_{B}\circ u_{\varphi_{B}}))\circ(H\otimes c_{H,B})\circ(((\mu_{H}\otimes\mu_{H})\circ(H\otimes c_{H,H}\otimes H)\circ(\delta_{H}\otimes\delta_{H}))\otimes B)\\&\footnotesize\textnormal{(by naturality of $c$)}\\=&\Phi'_{B}\circ(\mu_{H}\otimes B)\;\footnotesize\textnormal{(by the condition of coalgebra morphism for $\mu_{H}$ and definition of $\Phi'_{B}$)}.
\end{align*}
So we conclude that $(B,\Phi'_{B})$ is a left $H$-module. To finish the proof we have to compute that $\eta_{B}$ and $\mu_{B}$ are morphisms of left $H$-modules. On the one side,
\begin{align*}
&\Phi'_{B}\circ(H\otimes\eta_{B})\\=&u_{\varphi_{B}}\ast(\lambda_{B}\circ u_{\varphi_{B}})\;\footnotesize\textnormal{(by naturality of $c$ and definition of $u_{\varphi_{B}}$)}\\=&\varepsilon_{H}\otimes\eta_{B}\;\footnotesize\textnormal{(by \eqref{rmk2})}.
\end{align*}
On the other side,
\begin{align*} 
&\Phi'_{B}\circ(H\otimes\mu_{B})\\=&\mu_{B}\circ((\varphi_{B}\circ(H\otimes\mu_{B}))\otimes(\lambda_{B}\circ u_{\varphi_{B}}))\circ(H\otimes((B\otimes c_{H,B})\circ(c_{H,B}\otimes B)))\circ(\delta_{H}\otimes B\otimes B)\\&\footnotesize\textnormal{(by definition of $\Phi'_{B}$ and naturality of $c$)}\\=&\mu_{B}\circ((\mu_{B}\circ(\Phi'_{B}\otimes\varphi_{B})\circ(H\otimes c_{H,B}\otimes B)\circ(\delta_{H}\otimes B\otimes B))\otimes(\lambda_{B}\circ u_{\varphi_{B}}))\circ(H\otimes((B\otimes c_{H,B})\\&\circ(c_{H,B}\otimes B)))\circ(\delta_{H}\otimes B\otimes B)\;\footnotesize\textnormal{(by \eqref{Compat prima})}\\=&\mu_{B}\circ(\Phi'_{B}\otimes(\mu_{B}\circ (\varphi_{B}\otimes(\lambda_{B}\circ u_{\varphi_{B}}))\circ(H\otimes c_{H,B})))\circ(H\otimes((c_{H,B}\otimes H)\circ(H\otimes c_{H,B})\\&\circ(\delta_{H}\otimes B))\otimes B)\circ(\delta_{H}\otimes B\otimes B)\;\footnotesize\textnormal{(by associativity of $\mu_{B}$ and coassociativity of $\delta_{H}$)}\\=&\mu_{B}\circ(\Phi'_{B}\otimes\Phi'_{B})\circ(H\otimes c_{H,B}\otimes B)\circ(\delta_{H}\otimes B\otimes B)\;\footnotesize\textnormal{(by naturality of $c$ and definition of $\Phi'_{B}$)}.\qedhere
\end{align*}
\end{proof}
\begin{lemma}
Let's assume that $\sf{C}$ is a symmetric and let $(H,B,\varphi_{B})$ be a Hopf bracoid in $\sf{C}$ such that $\varphi_{B}$ is a coalgebra morphism. If $H$ is a cocommutative Hopf algebra, then $\Phi'_{B}$ is a coalgebra morphism.
\end{lemma}
\begin{proof}
It is straightforward to prove that $\varepsilon_{B}\circ\Phi'_{B}=\varepsilon_{H}\otimes\varepsilon_{B}$ due to the condition of coalgebra morphism for $\mu_{B}$ and $\varphi_{B}$, \eqref{u-antip} and \eqref{u1coalgeps}. Moreover,
\begin{align*}
&\delta_{B}\circ\Phi'_{B}\\=&(\mu_{B}\otimes\mu_{B})\circ(B\otimes c_{B,B}\otimes B)\circ((\delta_{B}\circ\varphi_{B})\otimes(\delta_{B}\circ\lambda_{B}\circ u_{\varphi_{B}}))\circ(H\otimes c_{H,B})\circ(\delta_{H}\otimes B)\;\footnotesize\textnormal{(by the condition}\\&\footnotesize\textnormal{of coalgebra morphism for $\mu_{B}$)}\\=&(\mu_{B}\otimes\mu_{B})\circ(B\otimes c_{B,B}\otimes B)\circ(((\varphi_{B}\otimes\varphi_{B})\circ(H\otimes c_{H,B}\otimes B)\circ(\delta_{H}\otimes\delta_{B}))\otimes(c_{B,B}\circ(\lambda_{B}\otimes\lambda_{B})\\&\circ(u_{\varphi_{B}}\otimes u_{\varphi_{B}})\circ\delta_{H}))\circ(H\otimes c_{H,B})\circ(\delta_{H}\otimes B)\;\footnotesize\textnormal{(by the condition of coalgebra morphism for $\varphi_{B}$, \eqref{a-antip} and \eqref{u1coalgdelta})}\\=&((\mu_{B}\circ (\varphi_{B}\otimes (\lambda_{B}\circ u_{\varphi_{B}})))\otimes(\mu_{B}\circ(\varphi_{B}\otimes(\lambda_{B}\circ u_{\varphi_{B}}))\circ(H\otimes c_{H,B})))\circ(H\otimes((B\otimes c_{H,H}\otimes H)\\&\circ(c_{H,B}\otimes c_{H,H})\circ(H\otimes c_{H,B}\otimes H)\circ(\delta_{H}\otimes B\otimes H))\otimes B)\circ(H\otimes H\otimes c_{H,B}\otimes B)\\&\circ(((H\otimes\delta_{H})\circ\delta_{H})\otimes\delta_{B})\;\footnotesize\textnormal{(by naturality of $c$, symmetry of $c$ and coassociativity of $\delta_{H}$)}\\=&((\mu_{B}\circ (\varphi_{B}\otimes (\lambda_{B}\circ u_{\varphi_{B}})))\otimes\Phi'_{B})\circ(H\otimes((c_{H,B}\otimes H)\circ(H\otimes c_{H,B})\circ((c_{H,H}\circ\delta_{H})\otimes B))\otimes B)\\&\circ(\delta_{H}\otimes\delta_{B})\;\footnotesize\textnormal{(by naturality of $c$ and definition of $\Phi'_{B}$)}\\=&((\mu_{B}\circ (\varphi_{B}\otimes (\lambda_{B}\circ u_{\varphi_{B}})))\otimes\Phi'_{B})\circ(H\otimes((c_{H,B}\otimes H)\circ(H\otimes c_{H,B})\circ(\delta_{H}\otimes B))\otimes B)\circ(\delta_{H}\otimes\delta_{B})\\&\footnotesize\textnormal{(by cocommutativity of $\delta_{H}$)}\\=&(\Phi'_{B}\otimes\Phi'_{B})\circ(H\otimes c_{H,B}\otimes B)\circ(\delta_{H}\otimes\delta_{B})\;\footnotesize\textnormal{(by coassociativity of $\delta_{H}$ and definition of $\Phi'_{B}$)}.
\qedhere
\end{align*}
\end{proof}
\begin{theorem} \label{monoidalHBrcd}
Let's assume that $\sf{C}$ is symmetric.
\emph{\sf{HBrcd}} is a symmetric monoidal category with base object $(K,K,id_{K})$, tensor functor defined as follows:
\begin{align*}
\otimes\colon \textnormal{\sf{HBrcd}}\times\textnormal{\sf{HBrcd}}&\longrightarrow \textnormal{\sf{HBrcd}}\\
((H,B,\varphi_{B}),(A,D,\varphi_{D}))&\longmapsto(H\otimes A,B\otimes D,\varphi_{B\otimes D}),
\end{align*}
where $\varphi_{B\otimes D}\coloneqq (\varphi_{B}\otimes\varphi_{D})\circ(H\otimes c_{A,B}\otimes D)$, and symmetry isomorphism
\begin{align*}
c_{(H,B,\varphi_{B}),(A,D,\varphi_{D})}\coloneqq(c_{H,A},c_{B,D})\colon(H\otimes A,B\otimes D,\varphi_{B\otimes D})\longrightarrow (A\otimes H,D\otimes B,\varphi_{D\otimes B}) 
\end{align*}
for all Hopf bracoids $(H,B,\varphi_{B})$ and $(A,D,\varphi_{D})$.
\end{theorem}
\begin{proof}
First of all, note that if $H=(H,\eta_{H},\mu_{H},\varepsilon_{H},\delta_{H},\lambda_{H})$ and $A=(A,\eta_{A},\mu_{A},\varepsilon_{A},\delta_{A},\lambda_{A})$ are Hopf algebras in the symmetric monoidal category $\sf{C}$, then it is straightforward to prove that $$H\otimes A=(H\otimes A,\eta_{H}\otimes\eta_{A},\mu_{H\otimes A},\varepsilon_{H}\otimes\varepsilon_{A},\delta_{H\otimes A},\lambda_{H}\otimes\lambda_{A})$$ is also a Hopf algebra in $\sf{C}$, where $\mu_{H\otimes A}$ and $\delta_{H\otimes A}$ are the usual tensor product and tensor coproduct, respectively. The same happens for $B\otimes D$, which is a Hopf algebra in $\sf{C}$ with antipode $\lambda_{B}\otimes\lambda_{D}$. 

Moreover, if $(B,\varphi_{B})$ is a left $H$-module and $(D,\varphi_{D})$ is a left $A$-module, then $$(B\otimes D,\varphi_{B\otimes D}\coloneqq (\varphi_{B}\otimes\varphi_{D})\circ(H\otimes c_{A,B}\otimes D))$$ is a left $H\otimes A$-module. Indeed:
\begin{itemize}
\item[(i)] It is clear that $\varphi_{B\otimes D}\circ(\eta_{H}\otimes\eta_{A}\otimes B\otimes D)=id_{B\otimes D}$ thanks to naturality of $c$ and identities $\varphi_{B}\circ(\eta_{H}\otimes B)=id_{B}$ and $\varphi_{D}\circ(\eta_{A}\otimes D)=id_{D}$.
\end{itemize}
\begin{itemize}
\item[(ii)] Let's see that $\varphi_{B\otimes D}\circ(H\otimes A\otimes\varphi_{B\otimes D})=\varphi_{B\otimes D}\circ(\mu_{H\otimes A}\otimes B\otimes D)$ also holds:
\begin{align*}
&\varphi_{B\otimes D}\circ(H\otimes A\otimes\varphi_{B\otimes D})\\=&((\varphi_{B}\circ(H\otimes\varphi_{B}))\otimes(\varphi_{D}\circ(A\otimes\varphi_{D})))\circ(H\otimes((H\otimes c_{A,B}\otimes A)\circ(c_{A,H}\otimes c_{A,B}))\otimes D)\\&\footnotesize\textnormal{(by naturality of $c$)}\\=&((\varphi_{B}\circ(\mu_{H}\otimes B))\otimes(\varphi_{D}\circ(\mu_{A}\otimes D)))\circ(H\otimes((H\otimes c_{A,B}\otimes A)\circ(c_{A,H}\otimes c_{A,B}))\otimes D)\;\footnotesize\textnormal{(by conditions}\\&\footnotesize\textnormal{of $H$-module and $A$-module)}\\=&\varphi_{B\otimes D}\circ(\mu_{H\otimes A}\otimes B\otimes D)\;\footnotesize\textnormal{(by naturality of $c$)}.
\end{align*}
\end{itemize}

So, to conclude that $(H\otimes A,B\otimes D,\varphi_{B\otimes D})$ is a Hopf bracoid, it is enough to see that \eqref{Eq. Bracoid compatibility} holds. To begin with, note that
\begin{equation}\label{Gamma tensor}
\Phi_{B\otimes D}=(\Phi_{B}\otimes\Phi_{D})\circ(H\otimes c_{A,B}\otimes D)
\end{equation}
as will be shown below. For the proof it will be necessary to take into account that
\begin{equation}\label{tensoru1}
u_{\varphi_{B\otimes D}}=u_{\varphi_{B}}\otimes u_{\varphi_{D}},
\end{equation}
which is trivial thanks to naturality of $c$. Then,
\begin{align*}
&\Phi_{B\otimes D}\\=&(\mu_{B}\otimes\mu_{D})\circ(B\otimes c_{D,B}\otimes D)\circ((\lambda_{B}\circ u_{\varphi_{B}})\otimes(\lambda_{D}\circ u_{\varphi_{D}})\otimes((\varphi_{B}\otimes\varphi_{D})\circ(H\otimes c_{A,B}\otimes D)))\\&\circ(((H\otimes c_{H,A}\otimes A)\circ(\delta_{H}\otimes\delta_{A}))\otimes B\otimes D)\;\footnotesize\textnormal{(by \eqref{tensoru1} and definitions of $\mu_{B\otimes D}$, $\varphi_{B\otimes D}$ and $\delta_{H\otimes A}$)}\\=&(\Phi_{B}\otimes(\mu_{D}\circ ((\lambda_{D}\circ u_{\varphi_{D}})\otimes\varphi_{D})))\circ(H\otimes((c_{A,B}\otimes A)\circ(A\otimes c_{A,B})\circ(\delta_{A}\otimes B))\otimes D)\;\footnotesize\textnormal{(by naturality  }\\&\footnotesize\textnormal{and symmetry of $c$ and definition of $\Phi_{B}$)}\\=&(\Phi_{B}\otimes\Phi_{D})\circ(H\otimes c_{A,B}\otimes D)\;\footnotesize\textnormal{(by naturality of $c$ and definition of $\Phi_{D}$)}.
\end{align*}
Therefore,
\begin{align*}
&\mu_{B\otimes D}\circ(\varphi_{B\otimes D}\otimes\Phi_{B\otimes D})\circ(H\otimes A\otimes c_{H\otimes A,B\otimes D}\otimes B\otimes D)\circ(\delta_{H\otimes A}\otimes B\otimes D\otimes B\otimes D)\\=&\mu_{B\otimes D}\circ(\varphi_{B\otimes D}\otimes((\Phi_{B}\otimes\Phi_{D})\circ(H\otimes c_{A,B}\otimes D)))\circ(H\otimes A\otimes c_{H\otimes A,B\otimes D}\otimes B\otimes D)\\&\circ(\delta_{H\otimes A}\otimes B\otimes D\otimes B\otimes D)\;\footnotesize\textnormal{(by \eqref{Gamma tensor})}\\=&((\mu_{B}\circ(\varphi_{B}\otimes\Phi_{B})\circ(H\otimes c_{H,B}\otimes B)\circ(\delta_{H}\otimes B\otimes B))\otimes(\mu_{D}\circ(\varphi_{D}\otimes\Phi_{D})\circ(A\otimes c_{A,D}\otimes D)\\&\circ(\delta_{A}\otimes D\otimes D)))\circ(H\otimes((B\otimes c_{A,B}\otimes D)\circ(c_{A,B}\otimes c_{D,B}))\otimes D)\;\footnotesize\textnormal{(by naturality of $c$ and $c$ symmetry)}\\=&((\varphi_{B}\circ(H\otimes\mu_{B}))\otimes(\varphi_{D}\circ(A\otimes\mu_{D}))\circ(H\otimes((B\otimes c_{A,B}\otimes D)\circ(c_{A,B}\otimes c_{D,B}))\otimes D)\;\footnotesize\textnormal{(by \eqref{Eq. Bracoid compatibility})}\\=&(\varphi_{B}\otimes\varphi_{D})\circ(H\otimes c_{A,B}\otimes D)\circ(H\otimes A\otimes ((\mu_{B}\otimes\mu_{D})\circ(B\otimes c_{D,B}\otimes D)))\;\footnotesize\textnormal{(by naturality of $c$)}\\=&\varphi_{B\otimes D}\circ(H\otimes A\otimes \mu_{B\otimes D}).
\end{align*}
So, we conclude that {\sf HBrcd} is monoidal. 

To prove that {\sf HBrcd} is also symmetric it is enough to see that $(c_{H,A},c_{B,D})$ is a morphism in the category {\sf HBrcd}. Indeed, thanks to naturality of $c$, $c_{H,A}$ and $c_{B,D}$ are Hopf algebra morphisms in $\sf{C}$ and, on the other hand, \eqref{Cond. Mor bracoid} holds because
\begin{align*}
&c_{B,D}\circ \varphi_{B\otimes D}\\=&(\varphi_{D}\otimes\varphi_{B})\circ(A\otimes c_{H,D}\otimes B)\circ(c_{H,A}\otimes c_{B,D})\;\footnotesize\textnormal{(by naturality of $c$ and $c$ symmetry)}\\=&\varphi_{D\otimes B}\circ(c_{H,A}\otimes c_{B,D}).\qedhere
\end{align*}
\end{proof}

\section{Hopf bracoids and 1-cocycles}\label{sect3}
This section will be devoted to analysing the existence of a relationship between 1-cocycles and Hopf bracoids. 1-cocycles are defined as follows:
\begin{definition}
Let $H,B$ be Hopf algebras in $\sf{C}$ and $\pi\colon H\rightarrow B$ a morphism in $\sf{C}$. Moreover, consider that $(B,\gamma_{B})$ is a left $H$-module algebra in $\sf{C}$. We will say that $\pi$ is a 1-cocycle in $\sf{C}$ if $\pi$ is a coalgebra morphism and the following condition holds:
\begin{equation}\label{CondCociclo}
\mu_{B}\circ(\pi\otimes\gamma_{B})\circ(\delta_{H}\otimes\pi)=\pi\circ\mu_{H}.
\end{equation}
\end{definition}
\begin{definition} Let $\pi\colon H\rightarrow B$ and $\tau\colon H'\rightarrow B'$ be 1-cocycles in $\sf{C}$. We will say that a pair $(f,g)$ is a morphism of 1-cocycles in $\sf{C}$ between $\pi$ and $\tau$ if $f\colon H\rightarrow H' $ and $g\colon B\rightarrow B'$ are Hopf algebra morphisms and the following conditions hold:
\begin{equation}\label{CondMorCoc1}
\tau\circ f=g\circ \pi,
\end{equation}
\begin{equation}\label{CondMorCoc2}
g\circ\gamma_{B}=\gamma_{B'}\circ(f\otimes g).
\end{equation}
\end{definition}
1-cocycles together with their respective morphisms give rise to a category that we will denote by $\sf{1C}$. If $\pi\colon H\rightarrow B$ is a 1-cocycle such that $\pi$ is an isomorphism, then we recover the notion of invertible 1-cocycle. Then, invertible 1-cocycles constitute a full subcategory of ${\sf 1C}$ that we will denote by ${\sf IC}$. Note that for every 1-cocycle $\pi\colon H\rightarrow B$ the equality 
\begin{equation}\label{etacocycle}
\pi\circ\eta_{H}=\eta_{B}
\end{equation}
always holds. The proof of this fact follows in exactly the same way as for invertible 1-cocycles (see \cite[Remark 1.11]{AGV}). 

Let's begin by building a functor $F$ from the category of 1-cocycles to Hopf bracoids.
\begin{theorem}\label{funtor F} There exists a functor $F\colon\textnormal{\sf{1C}}\longrightarrow\hbr$ defined on objects by $F\left(\pi\begin{array}{c}H\\\downarrow\\B\end{array}\right)=(H,B,\psi_{B})$, where $\psi_{B}\coloneqq \mu_{B}\circ(\pi\otimes\gamma_{B})\circ(\delta_{H}\otimes B)$, and on morphisms by the identity.
\end{theorem}
\begin{proof} 
First of all we have to see that $(B,\psi_{B})$ is a left $H$-module. On the one hand,
\begin{align*}
&\psi_{B}\circ(\eta_{H}\otimes B)\\=&\mu_{B}\circ((\pi\circ\eta_{H})\otimes(\gamma_{B}\circ(\eta_{H}\otimes B)))\;\footnotesize\textnormal{(by the condition of coalgebra morphism for $\eta_{H}$)}\\=&\mu_{B}\circ(\eta_{B}\otimes B)\;\footnotesize\textnormal{(by the conditions of left $H$-module and \eqref{etacocycle})}\\=&id_{B}\;\footnotesize\textnormal{(by unit property)}.
\end{align*}
On the other hand,
\begin{align*}
&\psi_{B}\circ(H\otimes\psi_{B})\\=&\mu_{B}\circ(\pi\otimes(\mu_{B}\circ(\gamma_{B}\otimes\gamma_{B})\circ(H\otimes c_{H,B}\otimes B)\circ(\delta_{H}\otimes B\otimes B)))\circ(\delta_{H}\otimes((\pi\otimes\gamma_{B})\circ(\delta_{H}\otimes B)))\\&\footnotesize\textnormal{(by the condition of morphism of left $H$-modules for $\mu_{B}$)}\\=&\mu_{B}\circ((\mu_{B}\circ(\pi\otimes\gamma_{B})\circ(\delta_{H}\otimes\pi))\otimes(\gamma_{B}\circ(\mu_{H}\otimes B)))\circ(((H\otimes c_{H,H}\otimes H)\circ(\delta_{H}\otimes\delta_{H}))\otimes B)\\&
\footnotesize\textnormal{(by naturality of $c$, associativity of $\mu_{B}$, coassociativity of $\delta_{H}$ and the conditions of left $H$-module)}\\=&\mu_{B}\circ(\pi\otimes\gamma_{B})\circ(((\mu_{H}\otimes\mu_{H})\circ(H\otimes c_{H,H}\otimes H)\circ(\delta_{H}\otimes\delta_{H}))\otimes B)\;\footnotesize\textnormal{(by \eqref{CondCociclo})}\\=&\psi_{B}\circ(\mu_{H}\otimes B)\;\footnotesize\textnormal{(by the condition of coalgebra morphism for $\mu_{H}$ and definition of $\psi_{B}$)}.
\end{align*}
So, let's see that \eqref{Eq. Bracoid compatibility} holds for $(H,B,\psi_{B})$. First note that, by the condition of morphism of left $H$-modules for $\eta_{B}$, it is straightforward to prove that
\begin{equation}\label{u1psiB}
u_{\psi_{B}}=\pi,
\end{equation}
and then, \begin{equation}\label{GammaB=gammaB}\Phi_{B}=\gamma_{B}.\end{equation} Indeed:
\begin{align*}
&\Phi_{B}\\=&\mu_{B}\circ((\lambda_{B}\circ \pi)\otimes(\mu_{B}\circ(\pi\otimes\gamma_{B})\circ(\delta_{H}\otimes B)))\circ(\delta_{H}\otimes B)\;\footnotesize\textnormal{(by definition of $\psi_{B}$ and \eqref{u1psiB})}\\=&\mu_{B}\circ((\mu_{B}\circ((\lambda_{B}\circ\pi)\otimes\pi)\circ\delta_{H})\otimes\gamma_{B})\circ(\delta_{H}\otimes B)\;\footnotesize\textnormal{(by coassociativity of $\delta_{H}$ and associativity of $\mu_{B}$)}\\=&\mu_{B}\circ(((\lambda_{B}\ast id_{B})\circ\pi)\otimes\gamma_{B})\circ(\delta_{H}\otimes B)\;\footnotesize\textnormal{(by the condition of coalgebra morphism for $\pi$)}\\=&\gamma_{B}\;\footnotesize\textnormal{(by \eqref{antipode}, the condition of coalgebra morphism for $\pi$ and (co)unit properties)}.
\end{align*}
As a consequence,
\begin{align*}
&\mu_{B}\circ(\psi_{B}\otimes\Phi_{B})\circ(H\otimes c_{H,B}\otimes B)\circ(\delta_{H}\otimes B\otimes B)\\=&\mu_{B}\circ((\mu_{B}\circ(\pi\otimes\gamma_{B})\circ(\delta_{H}\otimes B))\otimes\gamma_{B})\circ(H\otimes c_{H,B}\otimes B)\circ(\delta_{H}\otimes B\otimes B)\;\footnotesize\textnormal{(by \eqref{GammaB=gammaB} and definition of $\psi_{B}$)}\\=&\mu_{B}\circ(\pi\otimes(\mu_{B}\circ(\gamma_{B}\otimes\gamma_{B})\circ(H\otimes c_{H,B}\otimes B)\circ(\delta_{H}\otimes B\otimes B)))\circ(\delta_{H}\otimes B\otimes B)\;\footnotesize\textnormal{(by associativity of $\mu_{B}$}\\&\footnotesize\textnormal{and coassociativity of $\delta_{H}$)}\\=&\mu_{B}\circ(\pi\otimes\gamma_{B})\circ(\delta_{H}\otimes\mu_{B})\;\footnotesize\textnormal{(by the condition of morphism of left $H$-modules for $\mu_{B}$)}\\=&\psi_{B}\circ(H\otimes\mu_{B}).
\end{align*}
Therefore, we have that $F$ is well-defined on objects. Also, $F$ is well-defined on morphisms because if
$$(f,g)\colon\pi\begin{array}{c}H\\\downarrow\\B\end{array}\rightarrow\tau\begin{array}{c}H'\\\downarrow\\B'\end{array}$$ is a morphism in {\sf 1C}, then we have that $g\circ\psi_{B}=\psi_{B'}\circ(f\otimes g)$. Indeed:
\begin{align*}
&g\circ\psi_{B}\\=&\mu_{B'}\circ(g\otimes g)\circ(\pi\otimes\gamma_{B})\circ(\delta_{H}\otimes B)\;\footnotesize\textnormal{(by the condition of algebra morphism for $g$)}\\=&\mu_{B'}\circ(\tau\otimes\gamma_{B'})\circ(((f\otimes f)\circ\delta_{H})\otimes g)\;\footnotesize\textnormal{(by \eqref{CondMorCoc1} and \eqref{CondMorCoc2})}\\=&\psi_{B'}\circ(f\otimes g)\;\footnotesize\textnormal{(by the condition of coalgebra morphism for $f$ and definition of $\psi_{B'}$)}.\qedhere
\end{align*}
\end{proof}
Under the conditions of the previous theorem, we also obtain that $\psi_{B}$ is a coalgebra morphism if we require $\gamma_{B}$ to be a coalgebra morphism and that \eqref{gammaBclascoc} holds , i.e., $$(H,B,\psi_{B})=F\left(\pi\begin{array}{c}H\\\downarrow\\B\end{array}\right)$$ is an object in $\overline{\sf HBrcd}$.
\begin{lemma}\label{PsiBcoalg}
Let $\pi\colon H\rightarrow B$ be a 1-cocycle such that the equality
\begin{equation}\label{gammaBclascoc}
(\gamma_{B}\otimes H)\circ(H\otimes c_{H,B})\circ((c_{H,H}\circ\delta_{H})\otimes B)=(\gamma_{B}\otimes H)\circ(H\otimes c_{H,B})\circ(\delta_{H}\otimes B)
\end{equation}
holds. If $\gamma_{B}$ is a coalgebra morphism, then $\psi_{B}$ is a coalgebra morphism. 
\end{lemma}
\begin{proof}
First of all, note that it is straightforward to see that  $\varepsilon_{B}\circ\psi_{B}=\varepsilon_{H}\otimes\varepsilon_{B}$ because $\varepsilon_{B}\circ \gamma_{B}=\varepsilon_{H}\otimes\varepsilon_{B}$. On the other hand, we have that
\begin{align*}
&\delta_{B}\circ\psi_{B}\\=&(\mu_{B}\otimes\mu_{B})\circ(B\otimes c_{B,B}\otimes B)\circ((\delta_{B}\circ\pi)\otimes(\delta_{B}\circ\gamma_{B}))\circ(\delta_{H}\otimes B)\;\footnotesize\textnormal{(by the condition of coalgebra}\\&\footnotesize\textnormal{morphism for $\mu_{B}$)}\\=&(\mu_{B}\otimes\mu_{B})\circ(B\otimes c_{B,B}\otimes B)\circ (((\pi\otimes\pi)\circ\delta_{H})\otimes((\gamma_{B}\otimes\gamma_{B})\circ(H\otimes c_{H,B}\otimes B)\circ(\delta_{H}\otimes\delta_{B})))\\&\circ(\delta_{H}\otimes B)\;\footnotesize\textnormal{(by the condition of coalgebra morphism for $\pi$ and $\gamma_{B}$)}\\=&((\mu_{B}\circ(\pi\otimes B))\otimes(\mu_{B}\circ(\pi\otimes\gamma_{B})))\circ(H\otimes((\gamma_{B}\otimes H)\circ(H\otimes c_{H,B})\circ((c_{H,H}\circ\delta_{H})\otimes B))\otimes H\otimes B)\\&\circ(H\otimes H\otimes c_{H,B}\otimes B)\circ(((H\otimes\delta_{H})\circ\delta_{H})\otimes\delta_{B})\;\footnotesize\textnormal{(by naturality of $c$ and coassociativity of $\delta_{H}$)}\\=&((\mu_{B}\circ(\pi\otimes B))\otimes(\mu_{B}\circ(\pi\otimes\gamma_{B})))\circ(H\otimes((\gamma_{B}\otimes H)\circ(H\otimes c_{H,B})\circ(\delta_{H}\otimes B))\otimes H\otimes B)\\&\circ(H\otimes H\otimes c_{H,B}\otimes B)\circ(((H\otimes\delta_{H})\circ\delta_{H})\otimes\delta_{B})\;\footnotesize\textnormal{(by \eqref{gammaBclascoc})}\\=&(\psi_{B}\otimes(\mu_{B}\circ (\pi\otimes\gamma_{B})))\circ(H\otimes((c_{H,B}\otimes H)\circ(H\otimes c_{H,B})\circ(\delta_{H}\otimes B))\otimes B)\circ(\delta_{H}\otimes\delta_{B})\\&\footnotesize\textnormal{(by coassociativity of $\delta_{H}$ and definition of $\psi_{B}$)}\\=&(\psi_{B}\otimes\psi_{B})\circ(H\otimes c_{H,B}\otimes B)\circ(\delta_{H}\otimes\delta_{B})\;\footnotesize\textnormal{(by naturality of $c$ and definition of $\psi_{B}$)}.\qedhere
\end{align*}
\end{proof}
\begin{remark}\label{cc2} If ${\sf C}$ is symmetric, \eqref{gammaBclascoc} means that $(B,\gamma_{B})$ is in the cocommutativity class of $H$.
Note that if $\pi\colon H\rightarrow B$ is a 1-cocycle such that $H$ is cocommutative, then \eqref{gammaBclascoc} always holds.
\end{remark}
So, if we denote by $\overline{{\sf 1C}}$ to the full subcategory of {\sf 1C} whose objects are 1-cocycles $\pi\colon H\rightarrow B$ verifying that $\gamma_{B}$ is a coalgebra morphism, by $\overline{{\sf 1C}}^{\,\star}$ to the full subcategory of $\overline{{\sf 1C}}$ such that \eqref{gammaBclascoc} holds, and by $\overline{{\sf HBrcd}}^{\,\star}$ to the full subcategory of $\overline{{\sf HBrcd}}$ such that \eqref{clas coc} holds, then functor $F$ defined in Theorem \ref{funtor F} restricts to a functor
\[F'\colon \overline{{\sf 1C}}^{\,\star}\longrightarrow \overline{{\sf HBrcd}}^{\,\star}\]
because when $\pi\colon H\rightarrow B$ is a 1-cocycle in $\overline{{\sf 1C}}^{\,\star}$, then $\psi_{B}$ is a coalgebra morphism by Lemma \ref{PsiBcoalg}, and \eqref{clas coc} also holds thanks to \eqref{GammaB=gammaB} and \eqref{gammaBclascoc}.

Moreover, we are able to build a functor $G$ from Hopf bracoids to 1-cocycles.
\begin{theorem}\label{functorG}
There exists a functor $G\colon \overline{\textnormal{\sf{HBrcd}}}\longrightarrow \textnormal{\sf{1C}}$ acting on objects by $G((H,B,\varphi_{B}))=u_{\varphi_{B}}\begin{array}{c}H\\\downarrow\\B\end{array}$ and on morphisms by the identity.
\end{theorem}
\begin{proof}
Let $(H,B,\varphi_{B})$ be a Hopf bracoid in ${\sf C}$ such that $\varphi_{B}$ is a coalgebra morphism. Thanks to Theorem \ref{GammaB modulo algebra}, we know that $(B,\Phi_{B})$ is a left $H$-module algebra, and $u_{\varphi_{B}}$ is a coalgebra morphism by equations \eqref{u1coalgdelta} and \eqref{u1coalgeps} in Remark \ref{rmkcoalgu1}. So, it is enough to see that \eqref{CondCociclo} holds to conclude that $G$ is well-defined on objects. Indeed: 
\begin{align*}
&\mu_{B}\circ(u_{\varphi_{B}}\otimes\Phi_{B})\circ(\delta_{H}\otimes u_{\varphi_{B}})\\=&\varphi_{B}\circ(H\otimes u_{\varphi_{B}})\;\footnotesize\textnormal{(by \eqref{phiB despejar})}\\=&u_{\varphi_{B}}\circ\mu_{H}\;\footnotesize\textnormal{(by \eqref{phiBu1})}.
\end{align*}
In addition, let $(f,g)\colon (H,B,\varphi_{B})\rightarrow(H',B',\varphi_{B'})$ be a morphism in $\overline{{\sf HBrcd}}$. Let's see that $(f,g)$ is also a morphism in {\sf 1C} between the 1-cocycles $u_{\varphi_{B}}\begin{array}{c}H\\\downarrow\\B\end{array}$ and $u_{\varphi_{B'}}\begin{array}{c}H'\\\downarrow \\B'\end{array}$. On the one hand, we have to compute that \eqref{CondMorCoc1} holds. Indeed:
\begin{align}\label{proofGfunc}
&g\circ u_{\varphi_{B}}\\\nonumber=&\varphi_{B'}\circ(f\otimes (g\circ\eta_{B}))\;\footnotesize\textnormal{(by \eqref{Cond. Mor bracoid})}\\\nonumber=&\varphi_{B'}\circ(f\otimes\eta_{B'})\;\footnotesize\textnormal{(by the condition of algebra morphism for $g$)}\\\nonumber=&u_{\varphi_{B'}}\circ f.
\end{align}
On the other hand, we have to prove that \eqref{CondMorCoc2} is satisfied.
\begin{align*}
&g\circ\Phi_{B}\\=&\mu_{B'}\circ((g\circ\lambda_{B}\circ u_{\varphi_{B}})\otimes(g\circ\varphi_{B}))\circ(\delta_{H}\otimes B)\;\footnotesize\textnormal{(by the condition of algebra morphism for $g$)}\\=&\mu_{B'}\circ((\lambda_{B'}\circ u_{\varphi_{B'}}\circ f)\otimes(\varphi_{B'}\circ(f\otimes g)))\circ(\delta_{H}\otimes B)\;\footnotesize\textnormal{(by \eqref{morant}, \eqref{proofGfunc} and \eqref{Cond. Mor bracoid})}\\=&\Phi_{B'}\circ(f\otimes g)\;\footnotesize\textnormal{(by the condition of coalgebra morphism for $f$ and definition of $\Phi_{D}$)}.
\end{align*}
As a conclusion, $G$ is also well-defined on morphisms.
\end{proof}
Note that if $(H,B,\varphi_{B})$ is an object in $\overline{{\sf HBrcd}}^{\,\star}$, then $\Phi_{B}$ is a coalgebra morphism by Lemma \ref{GammaBcoalg} and it results clear that \eqref{gammaBclascoc} holds for the 1-cocycle $u_{\varphi_{B}}$ thanks to \eqref{clas coc}. Therefore, functor $G$ defined in previous theorem restricts to a functor 
\[G'\colon \overline{{\sf HBrcd}}^{\,\star}\longrightarrow \overline{{\sf 1C}}^{\,\star}.\]
\begin{remark}\label{commutativeFunctorDiags}
As we have just mentioned in the introduction, the categories ${\sf IC}$ and ${\sf HBr}$ are equivalent (see \cite[Theorem 2.7]{FGRR2}). This categorical equivalence is induced by the following pair of functors: On the one hand, if $\mathbb{H}$ is a Hopf brace, then it is well-known that $(H_{1},\Gamma_{H_{1}})$ is a left $H_{2}$-module algebra, so there exists a functor
\[{\sf E}\colon {\sf HBr}\longrightarrow {\sf IC}\]
defined on objects by ${\sf E}(\mathbb{H})=id_{H}\begin{array}{c}H_{2}\\\downarrow\\H_{1}\end{array}$ and on morphisms by ${\sf E}(x)=(x,x)$. On the other hand, given an invertible 1-cocycle $\pi\colon A\rightarrow H$, we can obtain a new Hopf algebra structure $H_{\pi}=(H,\eta_{H},\mu_{H}^{\pi},\varepsilon_{H},\delta_{H},\lambda_{H}^{\pi})$ whose product and antipode are defined as follows:
\[\mu_{H}^{\pi}\coloneqq \pi\circ\mu_{A}\circ(\pi^{-1}\otimes\pi^{-1}),\;\;\lambda_{H}^{\pi}\coloneqq \pi\circ\lambda_{A}\circ\pi^{-1}.\]
Therefore, there exists another functor 
\[{\sf Q}\colon {\sf IC}\longrightarrow {\sf HBr}\]
defined by ${\sf Q}\left(\pi\begin{array}{c}A\\\downarrow\\H\end{array}\right)=\mathbb{H}_{\pi}=(H,H_{\pi})$ on objects and on morphisms by ${\sf Q}((f,g))=g$.

Consider then the following diagrams of functors:

\begin{minipage}{0.4 \textwidth}
\[\xymatrix{&{\sf HBr}\ar[d]^-{{\sf T}}\ar[rr]^-{{\sf E}} & &{\sf IC}\ar[d]^-{{\sf j}}\\&\overline{{\sf HBrcd}}\ar[rr]^-{G} & &{\sf 1C}}\]
\end{minipage}
\begin{minipage}{0.6\textwidth}
\[\xymatrix{&{\sf IC}\ar[d]^-{{\sf j}}\ar[rr]^-{{\sf Q}} & &{\sf HBr}\ar[d]^-{{\sf T}}\\&{\sf 1C}\ar[rr]^-{F} & &{\sf HBrcd}}\]
\end{minipage}
where ${\sf j}$ denotes the inclusion functor of ${\sf IC}$ into ${\sf 1C}$. Also note that ${\sf T}(\mathbb{H})=(H_{2},H_{1},\mu_{H}^{2})$ is always an object in $\overline{{\sf HBrcd}}$ because $\mu_{H}^{2}$ is a coalgebra morphism. It results straightforward to prove that $G\circ{\sf T}={\sf j}\circ{\sf E}$. However, if $\pi\colon H\rightarrow B$ is an invertible 1-cocycle, then we have that 
\[(F\circ{\sf j})\left(\pi\begin{array}{c}H\\\downarrow\\B\end{array}\right)=F\left(\pi\begin{array}{c}H\\\downarrow\\B\end{array}\right)=(H,B,\psi_{B}),\]
and 
\[({\sf T}\circ{\sf Q})\left(\pi\begin{array}{c}H\\\downarrow\\B\end{array}\right)={\sf T}(\mathbb{B}_{\pi})=(B_{\pi},B,\mu_{B}^{\pi})\]
which are isomorphic objects in the category of ${\sf HBrcd}$. Indeed, consider $(\pi,id_{B})\colon (H,B,\psi_{B})\rightarrow (B_{\pi},B,\mu_{B}^{\pi})$: At first, note that $\pi\colon H\rightarrow B_{\pi}$ is a Hopf algebra morphism because
\begin{align*}
&\mu_{B}^{\pi}\circ(\pi\otimes\pi)=\pi\circ\mu_{H}\circ((\pi^{-1}\circ \pi)\otimes(\pi^{-1}\circ\pi))=\pi\circ\mu_{H},
\end{align*} 
and, moreover, \eqref{Cond. Mor bracoid} follows by
\begin{align*}
&id_{B}\circ \psi_{B}\\=&\mu_{B}\circ(\pi\otimes\gamma_{B})\circ(\delta_{H}\otimes (\pi\circ\pi^{-1}))\;\footnotesize\textnormal{(by $\pi$ isomorphism)}\\=&\pi\circ\mu_{H}\circ(H\otimes\pi^{-1})\;\footnotesize\textnormal{(by \eqref{CondCociclo})}\\=&\mu_{B}^{\pi}\circ(\pi\otimes id_{B})\;\footnotesize\textnormal{(by definition of $\mu_{B}^{\pi}$ and $\pi$ isomorphism)}.
\end{align*}
So, due to being $\pi$ and $id_{B}$ isomorphism in ${\sf C}$, $(\pi,id_{B})\colon (H,B,\psi_{B})\rightarrow (B_{\pi},B,\mu_{B}^{\pi})$ is an isomorphism in ${\sf HBrcd}$.
\end{remark}
\begin{theorem}\label{mainiso}
Categories $\overline{{\sf 1C}}^{\,\star}$ and $\overline{{\sf HBrcd}}^{\,\star}$ are isomorphic.
\end{theorem}
\begin{proof} At first we will prove that $G'\circ F'={\sf id}_{\overline{{\sf 1C}}^{\,\star}}$. Indeed, consider $\pi\colon H\rightarrow B$ an object in $\overline{{\sf 1C}}^{\,\star}$, we have that
\begin{align*}
&(G'\circ F')(\pi\begin{array}{c}H\\\downarrow\\ B\end{array})=G'((H,B,\psi_{B}))=u_{\psi_{B}}\begin{array}{c}H\\\downarrow\\B\end{array}=\pi\begin{array}{c}H\\\downarrow\\B\end{array}\;\footnotesize\textnormal{(by \eqref{u1psiB})},	
\end{align*}
and also $\Phi_{B}=\gamma_{B}$ by \eqref{GammaB=gammaB}. 

On the other hand, if $(H,B,\varphi_{B})$ is an object in $\overline{{\sf HBrcd}}^{\,\star}$, then
\begin{align*}
&(F'\circ G')((H,B,\varphi_{B}))=F'\left(u_{\varphi_{B}}\begin{array}{c}H\\\downarrow\\B\end{array}\right)=(H,B,\psi_{B}),
\end{align*}
where $\psi_{B}$ for the Hopf bracoid $F'\left(u_{\varphi_{B}}\begin{array}{c}H\\\downarrow\\B\end{array}\right)$ is given by
\begin{align*}
&\psi_{B}=\mu_{B}\circ(u_{\varphi_{B}}\otimes\Phi_{B})\circ(\delta_{H}\otimes B)=\varphi_{B}\;\footnotesize\textnormal{(by \eqref{phiB despejar})}.
\end{align*}
Then, $F'\circ G'={\sf id}_{\overline{{\sf HBrcd}}^{\,\star}}$.
\end{proof}
Denoting by ${\sf coc}\overline{{\sf HBrcd}}$ to the full subcategory of $\overline{{\sf HBrcd}}$ whose objects are Hopf bracoids $(H,B,\varphi_{B})$ with $H$ cocommutative, and by ${\sf coc}\overline{{\sf 1C}}$ to the full subcategory of $\overline{{\sf 1C}}$ whose objects are 1-cocycles $\pi\colon H\rightarrow B$ such that $H$ is cocommutative, the following result is a direct consequence of the previous theorem.
\begin{corollary}
Categories ${\sf coc}\overline{{\sf 1C}}$ and ${\sf coc}\overline{{\sf HBrcd}}$ are isomorphic.
\end{corollary}
\begin{proof}
It is a direct consequence of Theorem \ref{mainiso} taking into account Remarks \ref{cc1} and \ref{cc2}.
\end{proof}
\section*{Funding}
The  authors were supported by  Ministerio de Ciencia e Innovaci\'on of Spain. Agencia Estatal de Investigaci\'on. Uni\'on Europea - Fondo Europeo de Desarrollo Regional (FEDER). Grant PID2020-115155GB-I00: Homolog\'{\i}a, homotop\'{\i}a e invariantes categ\'oricos en grupos y \'algebras no asociativas.

Moreover, José Manuel Fernández Vilaboa and Brais Ramos Pérez were funded by Xunta de Galicia, grant ED431C 2023/31 (European FEDER support included, UE).

Also, Brais Ramos Pérez was financially supported by Xunta de Galicia Scholarship ED481A-2023-023.
\bibliographystyle{amsalpha}

\begin{thebibliography}{A}
\bibitem{CCH} J. N. Alonso Álvarez, J. M. Fernández Vilaboa and R. González Rodríguez, On the (co)-commutativity class of a Hopf algebra and crossed products in a braided category, \emph{Comm. Algebra} \textbf{29}(12) (2001) 5857-5878.
\bibitem{AGV} I. Angiono, C. Galindo and L. Vendramin, Hopf braces and Yang-Baxter operators, \emph{Proc. Am. Math. Soc.} \textbf{145}(5) (2017) 1981-1995.
\bibitem{BRZ} T. Brzezi\'nski, Trusses: Between braces and rings, \emph{Trans. Am. Math. Soc.} \textbf{372}(6) (2019) 4149-4176.
\bibitem{CCS} F. Catino, I. Colazzo and P. Stefanelli, Semi-braces and the Yang-Baxter equation, \emph{J. Algebra} \textbf{483} (2017) 163-187.
\bibitem{CMMS} F. Catino, M. Mazzotta, M. Miccoli and P. Stefanelli, Set-theoretic solutions of the Yang-Baxter equation associated to weak braces, \emph{Semigroup Forum} \textbf{104} (2022) 228-255.
\bibitem{FGRR2} J. M. Fernández Vilaboa, R. González Rodríguez, B. Ramos Pérez and A. B. Rodríguez Raposo, Modules for invertible 1-cocycles, Preprint (2023), arXiv:2311.05233.
\bibitem{RGON} R. González Rodríguez, The fundamental theorem of Hopf modules for Hopf braces, \emph{Linear Multilinear Algebra} \textbf{70}(20) (2022) 5146-5156.
\bibitem{GONROD} R. González Rodríguez and A. B. Rodríguez Raposo, Categorical equivalences for Hopf trusses and their modules, Preprint (2023), arXiv:2312.06520.
\bibitem{GV} L. Guarnieri and L. Vendramin, Skew braces and the Yang–Baxter equation, \emph{Math. Comput.} \textbf{86}(307) (2017) 2519-2534.
\bibitem{JS2}  A. Joyal and R. Street, Braided tensor categories, \emph{Adv. Math.} \textbf{102}(1) (1993) 20-78.
\bibitem{K}  C. Kassel, {\em Quantum Groups} (Springer-Verlag, New York, 1995).
\bibitem{Mac} S. Mac Lane, {\em Categories for the working mathematician} (Springer-Verlag, New York, 1998).
\bibitem{MLT} I. Martin-Lyons and P. J. Truman, Skew bracoids, \emph{J. Algebra}, \textbf{638} (2024) 751-787.
\bibitem{RADHA} D. E. Radford, \emph{Hopf algebras} (World Scientific, Singapore, 2012).
%\bibitem{Sch} P. Schauenburg, On the braiding on a Hopf algebra in a braided category, \emph{New York J. Math} \textbf{4} (1998) 259-263.
\end{thebibliography}

\end{document}